\documentclass[12pt,draftcls,onecolumn]{IEEEtran}
\usepackage{cite}

\ifCLASSINFOpdf
   \usepackage[pdftex]{graphicx}
   \graphicspath{{../pdf/}{../jpeg/}}
   \DeclareGraphicsExtensions{.pdf,.jpeg,.png}
\else
   \usepackage[dvips]{graphicx}
   \graphicspath{{../eps/}}
   \DeclareGraphicsExtensions{.eps}
\fi
\usepackage{amsmath}
\usepackage{algorithmic}
\usepackage{algorithm}

\usepackage{array}

\ifCLASSOPTIONcompsoc
  \usepackage[caption=false,font=normalsize,labelfont=sf,textfont=sf]{subfig}
\else
 \usepackage[caption=false,font=footnotesize]{subfig}
\fi
\usepackage{amssymb} 
\usepackage{epstopdf} 
\newtheorem{theorem}{Theorem}
\newtheorem{lemma}{Lemma}
\newtheorem{remark}{Remark}
\newtheorem{assumption}{Assumption}
\begin{document}
%
\title{Geometric Convergence for Distributed Optimization with Barzilai-Borwein Step Sizes}
%
%
%

\author{Juan~Gao,~ Xinwei~Liu,~Yu-Hong~Dai,~Yakui Huang,~and Peng Yang 
\thanks{This work is supported in part by the National Natural Science Foundation of China under Grants 11671116, 11271107 and 11701137, and in part by the Major Research Plan of the NSFC under Grant 91630202. (Corresponding author: Xinwei Liu.)}
\thanks{J. Gao is with the School of Artificial Intelligence, Hebei University of Technology, Tianjin 300401, China (e-mail: gaojuan514@163.com). X. W. Liu and Y. K. Huang are with the Institute of Mathematics, Hebei University of Technology, Tianjin 300401, China (e-mail: mathlxw@hebut.edu.cn; huangyakui2006@gmail.com). Y. H. Dai is with the LSEC, ICMSEC, Academy of Mathematics and Systems Science, Chinese Academy of Sciences, Beijing 100190, China, and also with the School of Mathematical Sciences, University of Chinese Academy of Sciences, Beijing 100049, China (e-mail: dyh@lsec.cc.ac.cn). P. Yang is with the State Key Laboratory of Reliability and Intelligence of Electrical Equipment and with the School of Artificial Intelligence, Hebei University of Technology, Tianjin 300401, China (e-mail: yangp@hebut.edu.cn).}}
\maketitle

\begin{abstract}
We consider a distributed multi-agent optimization problem over a time-invariant undirected graph, where each agent possesses a local objective function and all agents collaboratively minimize the average of all objective functions through local computations and communications among neighbors. Recently, a class of distributed gradient methods has been proposed that achieves both exact and geometric convergence when a constant step size is used. The geometric convergence of these methods is ensured for conservatively selected step sizes, but how to choose an appropriate step size while running the algorithms has not been fully addressed. The Barzilai-Borwein (BB) method is a simple and effective technique for step sizes and requires few storage and inexpensive computations. It has been widely applied in various areas. In this paper, we introduce the BB method to distributed optimization. Based on an adapt-then-combine variation of the dynamic average consensus approach and using multi-consensus inner loops, we propose a distributed gradient method with BB step sizes (DGM-BB-C). Our method computes the step size for each agent automatically which only depends on its local information and is independent of that for other agents, and the larger step sizes are always permissible. Our method can seek the exact optimum when the number of consensus steps stays constant. We prove that DGM-BB-C has geometric convergence to the optimal solution. Simulation results on a distributed sensing problem show that our method is superior to some advanced methods in terms of iterations, gradient evaluations, communications and the related cost framework. These results validate our theoretical discoveries.
\end{abstract}


%
\IEEEpeerreviewmaketitle

\section{Introduction}
%
%
%
%
\IEEEPARstart{D}{istributed} optimization among networked agents has recently received considerable attention due to its wide applications in various areas. Many practical problems can be described as distributed optimization which can be formulated as minimizing the average of all local objective functions $f(x)=\frac{1}{n}\sum^{n}_{i=1}f_{i}(x),$ where $f_{i}(x)$ is the local objective function only known by agent $i$, and $x\in \mathbb{R}^p$ is a decision vector. In particular, we focus on the scenarios where function $f_i$ is convex and continuously differentiable, and the agents are connected thought a time-invariant undirected graph. Typical examples include distributed machine learning  \cite{Cevher14}, \cite{Boyd11}, model predictive control \cite{Necoara08}, distributed spectrum sensing \cite{Bazerque10}, formation control \cite{Olshevsky10}, \cite{Ren06}, multi-agent target seeking \cite{Pu16}, wireless networks \cite{Cohen17}, \cite{Ling10}, power system control \cite{Gan13}, sensor networks \cite{Rabbat04}, and so on.

The research on distributed optimization can be traced back to the seminal works \cite{Bertsekas83}, \cite{Tsitsiklis86} in the 1980s. With the emergence of large-scale networks, distributed optimization methods have attracted tremendous interest of researchers \cite{Nedic09,JakoveticDNG14,Yuan16,Terelius11,MotaADMM13,ShiDADMM14,Mokhtari16,Eisen17}. The most popular choices are distributed gradient descent methods \cite{Nedic09,JakoveticDNG14,Yuan16}, distributed dual decomposition \cite{ Terelius11} and decentralized alternating direction method of multipliers (ADMM) \cite{MotaADMM13}, \cite{ShiDADMM14}. In \cite{ Terelius11,MotaADMM13,ShiDADMM14}, these methods are based on Lagrangian dual variables and have been shown to have nice convergence rates, but these methods need to minimize a Lagrangian-related function to obtain the primal and dual variables at each iteration, which causes a high computation burden. Compared with these methods, distributed gradient methods involve less computational cost and are easy to implement. In this paper, we focus our discussion on distributed gradient methods. For the earlier distributed gradient methods \cite{Nedic09,JakoveticDNG14,Yuan16}, in order to ensure to converge to the exact solution, these approaches need to utilize a diminishing step size, which may result in a slow rate of convergence. With a constant step size, these methods can be fast, but they only converge to a neighborhood of the exact solution. Recently, distributed gradient methods have achieved significant improvements, which provide us new variants that converge linearly to the exact solution, see \cite{Berahas18,Shi14,XiDEXTRA17,XuPhD16, Xu15,NedicDIGing17,Qu16,NedicATC17,Lu18,Xu18,Xin19} for example. These methods can be classified into three categories by the techniques used to ensure the exact convergence.
The first category of methods employs a multi-consensus inner loop strategy to the standard distributed gradient descent (DGD) method and the number of consensus steps is increasing at an appropriate rate \cite{Berahas18}. The second category of methods exploits the difference of two consecutive DGD iterates with different weight matrices to cancel the steady state error \cite{Shi14}, \cite{XiDEXTRA17}. The third category of methods replaces the local gradients with the tracking gradients built on dynamic consensus \cite{Zhu10}, which enables each agent to asymptotically learn the average of local gradients \cite{XuPhD16,Xu15,NedicDIGing17,Qu16,NedicATC17,Lu18,Xu18,Xin19}. By combining the gradient tracking technique with an adapt-then-combine strategy, some distributed gradient methods are capable of using uncoordinated constant step sizes for distributed optimization, such as \cite{XuPhD16,Xu15,NedicATC17,Lu18,Xu18}.

Note that these distributed gradient methods perform well in practice but their geometric convergence can only be ensured theoretically for conservatively selected step sizes. For example, EXTRA in \cite{Shi14} achieves the linear convergence rate with the step size being requested to be smaller than $\frac{2\mu}{L^2}$. In \cite{NedicATC17}, using the adapt-then-combine variant of DIGing, ATC-DIGing achieves faster convergence where the step sizes are restricted not to exceed $\frac{1}{2L}$. In \cite {Lu18}, the step sizes for the linear convergence rate in the strongly convex case should also not exceed $\frac{1}{2L}$. The step sizes of AsynDGM in \cite{Xu18} are assumed to be smaller than $\frac{1}{L}$. In \cite{Xin19}, the step sizes of FROST for the linear convergence rate in the strongly convex case are not greater than $\frac{1}{nL}$. In contrast, these methods can be faster practically with the larger step sizes as compared to their theoretical counterparts. It remains an open question to prove linear convergence under the larger step sizes. Moreover, how to choose an appropriate step size while running the algorithms has not been fully addressed in these distributed gradient methods. For these methods used in practice, a constant step size is usually tuned by hand, which may impact on their practical performance greatly.

Proposed by Barzilai and Borwein in \cite{BB88}, the BB method has rapidly emerged as a winning paradigm to attack nonlinear optimization problems. The BB method is a simple and effective technique for the choice of the step size and requires few storage and inexpensive computations. Moreover, the BB step size does not require any parameters and is automatically computed while running the algorithm. Many generalizations and variants of the BB method have largely been developed, some of them have provided the corresponding convergence analysis \cite{Raydan97,Dai05,Dai06,Dai13,Dai19}. Recent years have witnessed the successful applications of the BB method in image processing \cite{Wang07}, compressed sensing \cite{Wright09}, sparse reconstruction \cite{Wen10}, signal processing \cite{Liu11}, nonnegative matrix factorization \cite{Huang15} and machine learning \cite{Tan16}. These features and successful applications of the BB method motivate us to incorporate the BB method into distributed optimization. Note that \cite{Tan10} and \cite{Deroo12} use BB step sizes in developing their distributed gradient methods with a backtracking line search for solving different practical problems in distribution optimization, but they do not provide any theoretical justifications.

The aim of this paper is to present a distributed gradient method with the BB step sizes, which is capable of converging geometrically to the exact optimal solution. In addition, our convergence admits larger step sizes. Based on an adapt-then-combine variation of the dynamic average consensus approach and using multi-consensus inner loops, we propose a distributed gradient method with the BB step sizes (DGM-BB-C). The proposed method has geometric convergence to the exact optimal solution, which is generally much larger than the scalar $\frac{1}{2L}$ or $\frac{1}{L}$ requested by the aforementioned methods.

Our method computes the step size for each agent automatically which only depends on its local information and is independent of the step sizes for other agents, which is different from these works \cite{NedicATC17,Lu18,Xu18}, where the step sizes are obtained by perturbing the hand-optimal identical step size with random variables satisfying some distribution. Moreover, these methods in \cite{NedicATC17,Lu18,Xu18} converge linearly to the optimal solution under the assumptions that the heterogeneity of the step sizes is small enough (i.e., the step sizes are very close to each other), and the largest step size meets an upper bound which is a function on the heterogeneity. However, the above two assumptions are difficult to be both satisfied, because sufficiently small local step sizes may not ensure small heterogeneity. Our BB step sizes do not ensure small heterogeneity, so the geometric convergence of our method can not be proved by similar techniques. Without these assumptions, we prove that our method converges geometrically to the optimal solution. The proposed DGM-BB-C can seek the exact optimum when the number of consensus steps stays constant, which is different from NEAR-DGD$^{+}$ \cite{Berahas18}, where the method converges linearly to the exact optimum with a constant step size when the number of consensus steps is increasing at an appropriate rate.
We conduct numerical simulations on solving a distributed least squares problem. We show that DGM-BB-C is not sensitive to the choice of initial step sizes and using multi-consensus inner loops is helpful in improving the performance of DGM-BB-C. The simulation results show that DGM-BB-C has the lowest computation cost and communication cost to reach an $\epsilon$-optimal solution compared with some advanced methods.

This paper is organized as follows. Section II formally introduces the problem formulation along with the main assumptions and describes our novel algorithm DGM-BB-C. Section III states main convergence results. In Section IV, we provide some lemmas as the basis of the proof of convergence of DGM-BB-C. The proof of main results is presented in Section V. Section VI presents simulation results that show the efficiency of our algorithm. Finally, Section VII draws some conclusions.

\emph{Notations:}~Throughout the paper, all vectors are regarded as columns if not otherwise specified. Let each agent $i$ hold a local copy of the global variable $x^{i}\in \mathbb{R}^{p}$. Its value at iteration $k$ is denoted by $x_k^{i}$. We stack the raw version of all $x^i$ into a single matrix $\mathbf{x}=[x^1,x^2,...,x^n]^T\in \mathbb{R}^{n\times p}$. We introduce an aggregate objective function of local variables: $\mathbf{f}(\mathbf{x})=\sum^{n}_{i=1}f_i(x^i)$ where its gradient is denoted by $\nabla \mathbf{f}(\mathbf{x})=[\nabla f_1(x^1),\nabla f_2(x^2),...,\nabla f_n(x^n)]^T\in \mathbb{R}^{n\times p}$ . Let $\mathbf{1}$ refer to a column vector with all entries equal to one and $I$ represent the $n\times n$ identity matrix. For any matrix $\mathbf{s}\in \mathbb{R}^{n\times p}$, we denote its average vector as $\overline{s}=\frac{1}{n}\mathbf{1}^T\mathbf{s}\in \mathbb{R}^{1\times p}$ and its consensus violation as $\mathbf{s}-\mathbf{1}\overline{s}=\mathbf{s}-\frac{1}{n}\mathbf{1}\mathbf{1}^T\mathbf{s}$. In addition, the gradient of $f$ at $\overline{s}$ is denoted by $\nabla f(\overline{s})=\frac{1}{n}\sum^{n}_{i=1}\nabla f_i(\overline{s})\in\mathbb{R}^{1\times p}$. We use $\rho(\cdot)$ to denote the spectral radius of a square matrix. Let $\|\cdot\|$ represent standard Euclidean norm for vectors, and Frobenius norm for matrices.

\section{DGM-BB-C Development}

In this section, we firstly formulate the optimization problem along with the main assumptions. Then we present the BB step size and describe our DGM-BB-C method.

We consider a network of $n$ agents communicating over a time-invariant connected undirect graph $\mathcal{G}=(\mathcal{V,\mathcal{E}})$, where $\mathcal{V}$ is the set of agents, $\mathcal{E}$ is the collection of pairs, $(i,j),i,j\in \mathcal{V}$, such that agents $i$ and $j$ can exchange information with each other. For each agent $i$, let $\mathcal{N}_i=\{j|j\neq i, (i,j)\in \mathcal{E}\}$ be its set of neighbors. All the agents collaboratively solve the following distributed optimization problem:
 \begin{IEEEeqnarray}{rCl}\label{1}
      \min \limits_{x\in \mathbb{R}^p}f(x)=\frac{1}{n}\sum^{n}_{i=1}f_{i}(x),
      \end{IEEEeqnarray}
where each local objective function $f_i:\mathbb{R}^p \rightarrow \mathbb{R}$ is convex and differentiable, and known only by agent $i$.

For local objective functions $f_i$, $i=1,2,...,n$, we make the following two standard assumptions.
\begin{assumption}[Smoothness]
For each agent $i$, its local objective $f_i$ is differentiable and has Lipschitz continuous gradient, i.e., for any $x, y \in \mathbb{R}^p$, there exists a positive constant $L$ such that
  \begin{IEEEeqnarray}{rCl}\label{2}\nonumber
     \|\nabla f_i(x)-\nabla f_i(y)\|\leq L\|x-y\|.
      \end{IEEEeqnarray}
\end{assumption}
\begin{assumption}[Strong convexity]
For each agent $i$, its local objective $f_i$ is strongly convex, i.e., for any $x, y \in \mathbb{R}^p$, there exists a positive constant $\mu$ such that
      \begin{IEEEeqnarray}{rCl}\label{3} \nonumber
      f_i(x)\geq f_i(y)+ \nabla f_i(y)^T(x-y) +\frac{\mu}{2}\|x-y\|^2.
      \end{IEEEeqnarray}
\end{assumption}

Constants $L$ and $\mu$ satisfy $\mu\leq L$ (see \cite[Chapter 3]{Bubeck15} for details). It is immediately from Assumption 2 that problem $(\ref{1})$  has a unique optimal solution denoted by $x^*\in \mathbb{R}^p$.
\begin{assumption}
The graph $\mathcal{G}$ is connected and the weight matrix $W=[w_{ij}]\in \mathbb{R}^{n\times n}$ is doubly stochastic.
\end{assumption}

Let $\delta$ denote the spectral norm of the matrix $W-\frac{1}{n}\mathbf{1}\mathbf{1}^T$. By Assumption 3, we have $\delta<1$.
\subsection{The Barzilai-Borwein Step Size}
The iterative format of the original BB method for solving $(\ref{1})$ takes the following form:
 \begin{IEEEeqnarray}{rCl}\label{5}
     x_{k+1}=x_k-\alpha_k\nabla f(x_k),
      \end{IEEEeqnarray}
      in which $\alpha_k$ is computed by either
      \begin{IEEEeqnarray}{rCl}\label{6}
      \alpha_k=\frac{s_{k}^Ts_k}{s_{k}^Tz_k}
      \end{IEEEeqnarray}
      or
      \begin{IEEEeqnarray}{rCl}\label{7}
      \alpha_k=\frac{s_{k}^Tz_k}{z_{k}^Tz_k},
      \end{IEEEeqnarray}
where $s_k=x_k-x_{k-1}$ and $z_k=\nabla f(x_k)-\nabla f(x_{k-1})$ for $k\geq 1$. We now apply the BB step size to distributed optimization.
\subsection{DGM-BB-C Algorithm}
Notice that the step size cannot be straightly computed utilizing the formulae from $(\ref{6})$ and $(\ref{7})$, because distributed optimization methods never compute the average gradient $\nabla f(x_k)$. Therefore, we need to implement the BB method in distributed fashion. Briefly speaking, at each iteration $k$, we replace the global variables $x_k$ and $x_{k-1}$ in $(\ref{6})$ and $(\ref{7})$ with local variables $x^i_k$ and $x^i_{k-1}$. Correspondingly, we replace the average gradients $\nabla f(x_k)$ and $\nabla f(x_{k-1})$ by the local gradients $\nabla f(x^i_k)$ and $\nabla f(x^i_{k-1})$. In this way, we obtain the distributed BB step sizes:
    \begin{IEEEeqnarray}{rCl}\label{8}
       \alpha_{k}^i=\frac{ (s_{k}^i)^T s_{k}^i}{ (s_{k}^i)^T z_{k}^i}
        \end{IEEEeqnarray}
         or
          \begin{IEEEeqnarray}{rCl}\label{9}
         \alpha_{k}^i=\frac{(s_{k}^i)^T z_{k}^i}{(z_{k}^i)^T z_{k}^i},
       \end{IEEEeqnarray}
       where $ s_{k}^i=x_{k}^i-x_{k-1}^i, z^i_{k}=\nabla f_i(x_{k}^i)-\nabla f_i(x_{k-1}^i)$.

Now we introduce the distributed BB step sizes into distributed optimization. In our algorithm, at iteration $k$, each agent $i$ maintains three variables, namely, $x_k^i, y_k^i, \alpha_k^i\in \mathbb{R}^{p} $. Each agent carries out the two steps: local optimization step and dynamic average consensus step. For local optimization step, different from \cite{Nedic09}, we utilize the adapt-then-combine strategy for local optimization:
     \begin{IEEEeqnarray}{rCl}\label{10}
 x_{k+1}^i=\sum_{j\in \mathcal{N}_i\cup \{i\}}w_{ij}[x_k^j-\alpha_k^j y_k^j],
 \end{IEEEeqnarray}
 where $\alpha_k^i$ is computed by $(\ref{8})$ or $(\ref{9})$ and $y_k^i$ is an estimated gradient to be computed in the consensus step. For dynamic average consensus step, in order to ensure the algorithm with distributed BB step sizes to seek the exact optimal solution, we use an adapt-then-combine variation of the dynamic average consensus approach to track the average of the gradients of objective functions:
 \begin{IEEEeqnarray}{rCl}\label{11}
  y_{k+1}^i=\sum_{j\in \mathcal{N}_i\cup \{i\}}w_{ij}[y_{k}^j+\nabla f_j(x_{k+1}^j)-\nabla f_j(x_{k}^j)].
\end{IEEEeqnarray}
However, this approach does not work well because of the deviation of the two gradient estimates. Thus, we conduct multi-consensus inner loops to make sure estimated gradients are as close to the average gradient as possible. In particular, we use a multi-consensus inner loop strategy for local optimization step and dynamic average consensus step, respectively. Let $R$ be a positive integer which is the number of inner consensus iterations. We summarize the proposed DGM-BB-C in Algorithm 1.
\begin{algorithm}
\caption{DGM-BB-C for Undirect connected Graph}
\label{alg1}
 1:~\textbf{Initialization:}~for every agent $i\in  \mathcal{V}$, $x_0^i\in \mathbb{R}^p$, $y_0^i=\nabla f_i(x_0^i)$, $\alpha_0^i>0$.\\
 2:~\textbf{Local Optimization:}~for every agent $i\in  \mathcal{V}$, computes:
   $$x_{k+1}^i(0)=x_k^i-\alpha_k^i y_k^i,$$
     $$ x_{k+1}^i(r)=\sum_{j\in \mathcal{N}_i\cup \{i\}}w_{ij}x_{k+1}^j(r-1),~~r=1,2,...,R,$$
     where $\alpha_k^i$ is computed by $(\ref{8})$ or $(\ref{9})$, and set $ x_{k+1}^i= x_{k+1}^i(R)$.\\
 3:~\textbf{Dynamic Average Consensus:}~for every agent $i\in  \mathcal{V}$, computes:
       $$y_{k+1}^i(0)=y_{k}^i+\nabla f_i(x_{k+1}^i)-\nabla f_i(x_{k}^i),$$
     $$ y_{k+1}^i(r)=\sum_{j\in \mathcal{N}_i\cup \{i\}}w_{ij}y_{k+1}^j(r-1),~~r=1,2,...,R,$$
     and set $ y_{k+1}^i= y_{k+1}^i(R)$.\\
 4:~Set $k\rightarrow k+1$ and go to Step 2.
\end{algorithm}

Based on the previous notations, the DGM-BB-C method can be rewritten as the following compact matrix form:
\begin{IEEEeqnarray}{rCl}
       \mathbf{x}_{k+1}&=&W^{R}[\mathbf{x}_k-D_k\mathbf{y}_k],\label{13}\\
       \mathbf{y}_{k+1}&=&W^{R}[\mathbf{y}_k+\nabla \mathbf{f}(\mathbf{x}_{k+1})-\nabla \mathbf{f}(\mathbf{x}_k)],\label{14}
      \end{IEEEeqnarray}
where $D_k$ is a diagonal matrix and $[D_k]^{ii}=\alpha_k^i$ is computed by $(\ref{8})$ or $(\ref{9})$.
\begin{remark}
The number of inner consensus iterations $R$ does not need to be too large, just make sure that it satisfies some lower bound which depends on the objective function and the underlying network through $\mu, L, n, \delta$, and an adjustable parameter vector $c\in\mathbb{R}^3$, which will be discussed in detail in Section III.
\end{remark}
\begin{remark}
If we always set $R=1$ and $\alpha_k^i=\alpha^i$ ($\alpha^i$ is a constant with different values for different agent $i$) in DGM-BB-C instead of using $(\ref{8})$ or $(\ref{9})$, then it reduces to ATC-DIGing. From this point of view, our method can be considered as an extension of ATC-DIGing. There are two important differences between them.
\begin{enumerate}
  \item Our method computes the step size for each agent automatically which only depends on its local information and is independent of the step sizes for other agents. This is totally different from ATC-DIGing, where the step sizes are derived from perturbing the hand-optimal identical step size by random variables satisfying the uniform distribution and are required to be very close to each other.
  \item DGM-BB-C uses a multi-consensus inner loop strategy for local optimization step and dynamic average consensus step, respectively, whereas ATC-DIGing does not. Due to these modifications, the larger step sizes are always permissible in our method, i.e., the step sizes of DGM-BB-C can be allowed to be not less than $\frac{1}{L}$ and can reach $\frac{1}{\mu}$. In contrast, the step sizes of ATC-DIGing can not be greater than $\frac{1}{2L}$ theoretically. These differences ensure DGM-BB-C has a better numerical performance than ATC-DIGing.
\end{enumerate}
\end{remark}
\section{Main Results}
In this section, we present main convergence results of DGM-BB-C. We start the convergence analysis by deriving the range of the distributed BB step sizes in DGM-BB-C, and then provide several important lemmas. Finally, based on these lemmas, Theorem 1 presents the geometric convergence of DGM-BB-C.
\begin{lemma}
Under Assumptions 1-2, for all $ k\geq 0$ and every agent $i \in \mathcal{V}$, the BB step size $\alpha_k^i$ computed by $(\ref{8})$ or $(\ref{9})$ in DGM-BB-C satisfies
\begin{IEEEeqnarray}{rCl}
 \frac{1}{L}\leq\alpha_k^i \leq \frac{1}{\mu}.
\end{IEEEeqnarray}
\end{lemma}
\begin{IEEEproof}
Firstly, we give the proof of bounds on the BB step size $\alpha_k^i$ computed by $(\ref{8})$.
By the strong convexity of local objective function $f_i$, we obtain
  \begin{IEEEeqnarray}{rCl}
(x_{k}^{i}-x_{k-1}^{i})^{T}(\nabla f_{i}(x_{k}^{i})-\nabla f_{i}(x_{k-1}^{i}))\geq \mu\|x_{k}^{i}-x_{k-1}^{i}\|^{2}.\IEEEnonumber
\end{IEEEeqnarray}
  Thus, we have the upper bound for each BB step size $\alpha_k^i$ computed by $(\ref{8})$ since
  \begin{IEEEeqnarray}{rCl}\label{15}
\alpha_k^i&=&\frac{(x_k^i-x_{k-1}^i)^T(x_k^i-x_{k-1}^i)}{(x_k^i-x_{k-1}^i)^T(\nabla f_i(x_k^i)-\nabla f_i(x_{k-1}^i))}\nonumber\\
 &\leq& \frac{\|x_k^i-x_{k-1}^i\|^2}{\mu\|x_k^i-x_{k-1}^i\|^2}\nonumber\\
 &=&\frac{1}{\mu}.
\end{IEEEeqnarray}
By the Cauchy inequality and the $L$-Lipschitz continuity of $\nabla f_i(x)$, one can derive that $\alpha_k^i$ is uniformly lower bounded due to
\begin{IEEEeqnarray}{rCl}
\alpha_k^i&=&\frac{(x_k^i-x_{k-1}^i)^T(x_k^i-x_{k-1}^i)}{(x_k^i-x_{k-1}^i)^T(\nabla f_i(x_k^i)-\nabla f_i(x_{k-1}^i))}\nonumber\\
&\geq& \frac{\|x_k^i-x_{k-1}^i\|^2}{\|x_k^i-x_{k-1}^i\|\|\nabla f_i(x_k^i)-\nabla f_i(x_{k-1}^i)\|}\nonumber\\
&\geq& \frac{\|x_k^i-x_{k-1}^i\|^2}{L\|x_k^i-x_{k-1}^i\|^2}\nonumber\\
&=&\frac{1}{L}.
\end{IEEEeqnarray}

Now, we give the proof of bounds on the BB step size $\alpha_k^i$ computed by $(\ref{9})$.
By the $L$-Lipschitz continuity of $\nabla f_i(x)$, we obtain
  \begin{IEEEeqnarray}{rCl}
(x_k^i-x_{k-1}^i)^T(\nabla f_i(x_k^i)-\nabla f_i(x_{k-1}^i))
\geq \frac{1}{L}\|\nabla f_i(x_k^i)-\nabla f_i(x_{k-1}^i)\|^2.
\end{IEEEeqnarray}
 Thus, $\alpha_k^i$ computed by $(\ref{9})$ is uniformly lower bounded by $\frac{1}{L}$. In terms of the strong convexity of local objective function $f_i$, we have $(s^i_k)^Tz^i_k>0$ . Hence, $\alpha_k^i$ computed by $(\ref{9})$ will not be greater than $\frac{1}{\mu}$. The desired result then follows.
\end{IEEEproof}

 We define $\alpha_{\max}\!\!=\!\max\limits_{k\geq 0}\{\alpha_k^i\}$ and $\overline{\alpha}_{\max}\!\!=\!\max\limits_{k\geq 0}\left\{\frac{1}{n}\!\sum^{n}_{i=1}\alpha_k^{i}\right\}$. It is immediately from Lemma 1 that
  \begin{IEEEeqnarray}{rCl}\label{16}
\frac{1}{L}\leq \alpha_{\max}, \overline{\alpha}_{\max}\leq \frac{1}{\mu}.
\end{IEEEeqnarray}
  By the above proof, we can see that $\alpha_k^i$ computed by $(\ref{8})$ is a longer step size while $\alpha_k^i$ computed by $(\ref{9})$ is a shorter one.

In the next lemma, we establish bounds on $\|\mathbf{x}_{k+1}-\mathbf{1}\overline{x}_{k+1}\|, \|\mathbf{y}_{k+1}-\mathbf{1}\overline{y}_{k+1}\|$ and $\|\overline{x}_{k+1}-(x^*)^T\|$ in terms of linear combinations of their past values.
\begin{lemma}
If $\frac{1}{n}\sum^{n}_{i=1}\alpha_{k}^{i}<\frac{2}{L}$ for all $k$, the following linear time invariant system inequality holds:
 \begin{IEEEeqnarray}{rCl}\label{17}
\mathbf{v}_{k+1}\leq G_{k}^{\alpha} \mathbf{v}_{k},~ \forall k
\end{IEEEeqnarray}
where $\mathbf{v}_{k}\in \mathbb{R}^{3},G_{k}^{\alpha}\in \mathbb{R}^{3\times 3}$  are defined as
\begin{equation}\label{23}\nonumber
 \mathbf{v}_{k}= \left[\begin{IEEEeqnarraybox*}[][c]{,c,}
    \|\mathbf{x}_{k}-\mathbf{1}\overline{x}_k\|    \\
         \|\mathbf{y}_{k}-\mathbf{1}\overline{y}_k\|   \\
          \|\overline{x}_{k}-(x^*)^T\|
\end{IEEEeqnarraybox*}\right],
\end{equation}
\begin{equation}
 G^{\alpha}_{k}\!=\! \left[\!\begin{IEEEeqnarraybox*}[][c]{,c/c/c,}
   \delta^R\!+\!\delta^RL\alpha_{\max}   & \delta^R\alpha_{\max}  &  \delta^R L \sqrt{n}\alpha_{\max}  \\
         2\delta^R L\!+\!\delta^RL^2\alpha_{\max} & \delta^R\!+\!\delta^RL\alpha_{\max}  &\delta^RL^2\sqrt{n}\alpha_{\max} \\
         \frac{L}{\sqrt{n}}\alpha_{\max}   & \frac{1}{\sqrt{n}}\alpha_{\max}  & \lambda
\end{IEEEeqnarraybox*}\!\right],\nonumber
\end{equation}
 $\lambda=\max\{|1-\frac{\mu}{n}\sum^{n}_{i=1}\alpha_{k}^{i}|, |1-\frac{L}{n}\sum^{n}_{i=1}\alpha_{k}^{i}|\}$, $R$ is the number of inner consensus iterations, and $\delta$ is the spectral norm of the matrix $W-\frac{1}{n}\mathbf{1}\mathbf{1}^T$.
 \end{lemma}
 \begin{IEEEproof}
 See Section V.
 \end{IEEEproof}

 Note that a linear iterative relation between $\mathbf{v}_{k+1}$ and $\mathbf{v}_{k}$ with matrix $G_{k}^{\alpha}$ is established in $(\ref{17})$. From Lemma 1, it follows that $\frac{1}{n}\sum^{n}_{i=1}\alpha_{k}^{i}\geq \frac{1}{L}$. If $\frac{1}{n}\sum^{n}_{i=1}\alpha_{k}^{i}\leq \frac{2}{L}-\frac{\mu}{L^2}$, then $\lambda\leq 1-\frac{\mu}{L}$. Hence, for all $k$, we have $G_k^{\alpha}\preceq G^{\alpha}$ with
\begin{equation}\label{24}
 G^{\alpha}\!=\! \left[\!\begin{IEEEeqnarraybox*}[][c]{,c/c/c,}
   \delta^R\!+\!\delta^RL\alpha_{\max}   & \delta^R\alpha_{\max}  &  \delta^R L \sqrt{n}\alpha_{\max}  \\
         2\delta^R L\!+\!\delta^RL^2\alpha_{\max} & \delta^R\!+\!\delta^RL\alpha_{\max}  &\delta^RL^2\sqrt{n}\alpha_{\max} \\
         \frac{L}{\sqrt{n}}\alpha_{\max}   & \frac{1}{\sqrt{n}}\alpha_{\max}  & 1-\frac{\mu}{L}
\end{IEEEeqnarraybox*}\!\right],
\end{equation}
where $\preceq$ denotes entry-wise less than or equal to. Thus, we get that $\rho(G_k^{\alpha})\leq \rho(G^{\alpha})$ \cite [Theorem 8.1.18]{Horn12}. If $\rho(G^{\alpha})<1$, then $(G^{\alpha})^k$ converges linearly to $0$ at rate $O(\rho(G^{\alpha})^k)$ \cite [Theorem 5.6.12]{Horn12}, in which case $\|\mathbf{v}_k\|$ also converges linearly to $0$ at rate $O(\rho(G^{\alpha})^k)$. In order to analyze the convergence rate of DGM-BB-C, we need the following lemma, which provides conditions that ensure  $\rho(G^{\alpha})<1$.
 \begin{lemma}
 (\!\!\cite[Corollary 8.1.29]{Horn12}) Let $A\in\mathbb{R}^{n\times n}$ be a nonnegative matrix and $\omega\in\mathbb{R}^{n}$ be a positive vector. If $A\omega <q\omega$, then $\rho(A)<q$.
\end{lemma}

 We now show that when the largest step size $\alpha_{\max}$ satisfies $(\ref{16})$, with the appropriate lower bound on the number of inner consensus iterations $R$, the spectral radius of $G^{\alpha}$ is less than $1$.
\begin{lemma}
Suppose that the Assumptions 1-3 hold. Consider the matrix $G^{\alpha}$ defined in $(\ref{24})$ with the largest step size $\alpha_{\max}$ satisfying $(\ref{16})$. Let $C=Lc_1+c_2+L\sqrt{n}c_3$ and \\$\Delta \equiv \min\{\frac{\mu c_1}{C+\mu c_1}, \frac{\mu c_2}{LC+2L\mu c_1+\mu c_2}\}$, where $c_1, c_2$ are any positive scalars and $c_3>\frac{L(Lc_1+c_2)}{\mu^2\sqrt{n}}$. If $R\geq \lceil\frac{ \ln \Delta}{\ln \delta}\rceil+\varphi(\frac{ \ln \Delta}{\ln \delta})$ where
$\varphi(x)=\left\{\begin{IEEEeqnarraybox}[\relax][c]{ls}
1, \:\: &for $x \in \mathbb{N}_{+}$\\
0,  &for $x \notin \mathbb{N}_{+}$%
\end{IEEEeqnarraybox}\right.$, then $\rho(G^{\alpha})<1$.
\end{lemma}
\begin{IEEEproof}
In light of Lemma 3, we derive the lower bound on the number of inner consensus iterations $R$ and a positive vector $c=[c_1,c_2,c_3]^T$ from
\begin{equation}\label{27}
 G^{\alpha}\left[\begin{IEEEeqnarraybox*}[][c]{,c,}
         c_1\\
         c_2\\
         c_3
\end{IEEEeqnarraybox*}\!\right]<\left[\begin{IEEEeqnarraybox*}[][c]{,c,}
         c_1\\
         c_2\\
         c_3
\end{IEEEeqnarraybox*}\!\right],
\end{equation}
which is equivalent to the following set of inequalities
\begin{eqnarray}\label{28}
\left\{\!\begin{array}{c}
         (\delta^RLc_1\!+\!\delta^Rc_2\!+\!\delta^RL \sqrt{n}c_3)\alpha_{\max} \!<\! c_1(1\!-\!\delta^R),   \\
         \!\delta^RL(Lc_1\!+\!c_2\!+\!L\sqrt{n}c_3)\alpha_{\max}\!<\!(1\!-\!\delta^R)c_2\!-\!2\delta^RLc_1, \\
         (\frac{Lc_1}{\sqrt{n}}\!+\!\frac{c_2}{\sqrt{n}})\alpha_{\max}\!<\!\frac{c_3\mu}{L}.\\
                       \end{array}\right.
\end{eqnarray}
Since the right hand side of the second inequality in $(\ref{28})$ has to be positive, we obtain
\begin{IEEEeqnarray}{rCl}\label{29}
\delta^R<\frac{c_2}{2Lc_1+c_2}.
\end{IEEEeqnarray}
 It follows from $(\ref{28})$ that
\begin{IEEEeqnarray}{rCl}\label{30}
\alpha_{\max}<\widehat{\alpha},
\end{IEEEeqnarray}
where
\begin{IEEEeqnarray}{rCl}\label{31}
\widehat{\alpha}\equiv\min\left\{\frac{(1\!-\!\delta^R)c_1}{\delta^RC},\frac{(1\!-\!\delta^R)c_2\!-\!2\delta^RLc_1}{\delta^RLC},\frac{\mu\sqrt{n}c_3}{L(L c_1\!+\!c_2)}\right\},\nonumber\\
\end{IEEEeqnarray}
 $C=Lc_1+c_2+L\sqrt{n}c_3$, $c_1, c_2, c_3$ are positive constants and the range of $\delta^R$ is given in $(\ref{29})$.
 Since the largest step size $\alpha_{\max}$ satisfies $(\ref{16})$, one has
  \begin{IEEEeqnarray}{rCl}\label{32}
\frac{1}{L}\leq \alpha_{\max}\leq \frac{1}{\mu}.
\end{IEEEeqnarray}
 In order to ensure that the range of the largest step size given in $(\ref{32})$ is contained in the range given in $(\ref{30})$, we require $\widehat{\alpha}>\frac{1}{\mu}$. That is,
\begin{IEEEeqnarray}{rCl}\label{33}
  \min\! \left\{\!\frac{(1\!-\!\delta^R)c_1}{\delta^RC},\frac{(1\!-\!\delta^R)c_2\!-\!2\delta^RLc_1}{\delta^RLC},\frac{\mu\sqrt{n}c_3}{L(L c_1\!+\!c_2)}\!\right\}\!>\!\frac{1}{\mu}.\nonumber\\
\end{IEEEeqnarray}
Combining $(\ref{29})$ with $(\ref{33})$, we get
\begin{eqnarray}\label{320}
\delta^R <\min\left\{\frac{\mu c_1}{C+\mu c_1},\frac{\mu c_2}{LC+2L\mu c_1+\mu c_2}\right\} \equiv\Delta,
\end{eqnarray}
and
\begin{eqnarray} \label{321}
c_3>\frac{L(Lc_1+c_2)}{\mu^2\sqrt{n}}.
\end{eqnarray}
It follows from $(\ref{320})$ that $R\geq \lceil\frac{ \ln \Delta}{\ln \delta}\rceil+\varphi(\frac{ \ln \Delta}{\ln \delta})$ where
$\varphi(x)=\left\{\begin{IEEEeqnarraybox}[\relax][c]{ls}
1, \:\: &for $x \in \mathbb{N}_{+}$\\
0,  &for $x \notin \mathbb{N}_{+}$%
\end{IEEEeqnarraybox}\right.$.
The desired result follows.
\end{IEEEproof}
\begin{remark}
We observe that the parameters $\mu, L, n, \delta$ can be determined when the problem and the underlying network under study are given. With these fixed parameters, the upper bound on the largest step size in $(\ref{31})$ depends on what the number of inner consensus iterations $R$ and the adjustable parameter vector $c$ are chosen. In order to ensure that the range of the largest step size given in $(\ref{16})$ is contained in the range given in $(\ref{30})$, we notice that the values of the first and second terms in the minimum of $(\ref{31})$ monotonically increase as $R$ grows. As a result, we can first obtain a lower bound on the number of inner consensus iterations $R$ relying on parameter $c$ to make sure that these terms are greater than $\frac{1}{\mu}$. Then we pick the suitable parameter $c$, which has numerous options. One of the feasible methods is to obtain the parameter $c$ by minimizing the lower bound on the number of inner consensus iterations $R$ under satisfying $(\ref{321})$. As long as guaranteeing $R$ to satisfy the lower bound and selecting the parameter $c$ to satisfy $(\ref{321})$, we can always make sure that $\widehat{\alpha}$ is greater than the upper bound on the largest step size $\frac{1}{\mu}$ given in $(\ref{16})$. The numerical experiments in Section VI also show that the number of inner consensus iterations $R$ does not need to be large enough to guarantee that the range of the largest step size $\alpha_{\max}$ given in $(\ref{16})$ is contained in the range given in $(\ref{30})$.
\end{remark}
\begin{remark}
Although the proof of Lemma 4 is similar to those used in some exiting distributed optimization methods such as \cite{Qu16,Xin19,XiADD18}, it is essentially different between them. With the range of the largest step size given, we derive the lower bound on the number of inner consensus iterations $R$ such that the spectral radius of $G^{\alpha}$ is less than $1$, whereas in these methods \cite{Qu16,Xin19,XiADD18}, the authors derive the range of the largest step size or the step size such that its counterpart is also satisfied.
\end{remark}

We now present the main convergence result of this paper in Theorem 1,
which shows the geometric convergence rate of DGM-BB-C to the optimal solution.
\begin{theorem}
Let the Assumptions 1-3 hold. If $\frac{1}{n}\sum^{n}_{i=1}\alpha_{k}^{i}\leq\frac{2}{L}-\frac{\mu}{L^2}$ for all $k$, and $\alpha_{\max}$, $R$ are used as that in Lemma 4, the sequence $\{\mathbf{x}_k\}$ generated by DGM-BB-C converges exactly to the unique optimizer $\mathbf{x}^*=\mathbf{1}(x^*)^T$ at a geometric rate, i.e., there exists some positive constant $M>0$ such that, for any $k$
\begin{IEEEeqnarray}{rCl}
\|\mathbf{x}_k-\mathbf{x}^*\|\leq M(\rho(G^{\alpha})+\xi)^k,
\end{IEEEeqnarray}
where $\xi$ is a arbitrarily small constant.
\end{theorem}
\begin{IEEEproof}
Applying $(\ref{17})$ recursively, we get
\begin{IEEEeqnarray}{rCl}\label{34}
\mathbf{v}_k \leq \left(\prod^{k-1}_{l=0}G^{\alpha}_{l}\right)\mathbf{v}_0,
\end{IEEEeqnarray}
Taking the norm on both sides of the above relation gives
 \begin{IEEEeqnarray}{rCl}\label{35}
\|\mathbf{v}_k\| &\leq& \left(\prod^{k-1}_{l=0}\rho(G^{\alpha}_l)\right)\|\mathbf{v}_0\|,\IEEEnonumber\\
&\leq& \rho(G^{\alpha})^k\|\mathbf{v}_0\|,
\end{IEEEeqnarray}
Denote $v_k=\sum^{k-1}_{l=0}\|\mathbf{v}_l\|$, then $(\ref{35})$ can be written as
 \begin{IEEEeqnarray}{rCl}\label{36}
\|\mathbf{v}_k\|=v_{k+1}-v_k \leq \rho(G^{\alpha})^{k}\|\mathbf{v}_0\|,
\end{IEEEeqnarray}
 which implies that $v_{k+1}\leq v_k + \rho(G^{\alpha})^{k}\|\mathbf{v}_0\|.$ By Lemma 4, one has that $v_k$ converges and therefore is bounded. For any $\gamma \in (\rho(G^{\alpha}),1)$, it follows from $(\ref{36})$ that
  \begin{IEEEeqnarray}{rCl}\label{37}
\lim_{k\rightarrow \infty}\frac{\|\mathbf{v}_k\|}{\gamma^k}\leq \lim_{k\rightarrow \infty}\frac{\rho(G^{\alpha})^{k}\|\mathbf{v}_0\|}{\gamma^k} \leq \|\mathbf{v}_0\|.
\end{IEEEeqnarray}
Therefore, $\|\mathbf{v}_k\|=O(\gamma^k)$. That is, there exists some positive constant $T$ such that, for all $k$
 \begin{IEEEeqnarray}{rCl}\label{38}
\|\mathbf{v}_k\|\leq T(\rho(G^{\alpha})+\xi)^k,
\end{IEEEeqnarray}
where $\xi$ is a arbitrarily small constant. Moreover, since
\begin{IEEEeqnarray}{rCl}\label{39}
\|\mathbf{x}_{k}-\mathbf{x}^{*}\|&\leq&\|\mathbf{x}_{k}-\mathbf{1}\overline{x}_{k}\|+\|\mathbf{1}\overline{x}_{k}-\mathbf{x}^{*}\|\IEEEnonumber\\
&\leq& \|\mathbf{x}_{k}-\mathbf{1}\overline{x}_{k}\|+ \sqrt{n}\|\overline{x}_{k}-(x^{*})^T\|\IEEEnonumber\\
&\leq& (1+\sqrt{n})\|\mathbf{v}_{k}\|,
\end{IEEEeqnarray}
by combining $(\ref{38})$ with $(\ref{39})$, we obtain that
\begin{IEEEeqnarray}{rCl}\label{40}
\|\mathbf{x}_k-\mathbf{x}^*\|\leq (1+\sqrt{n})T(\rho(G^{\alpha})+\xi)^k.
\end{IEEEeqnarray}
 The desired result follows immediately by letting $M=(1+\sqrt{n})T$.
\end{IEEEproof}
\begin{remark}
Theorem 1 shows that the sequence $\{\mathbf{x}_k\}$ converges geometrically to the optimal solution $\mathbf{x}^*$ with proper choices of $R$ and $c$. DGM-BB-C can seek the exact optimum theoretically in which the number of consensus steps does not need to be increasing at an appropriate rate. We will show the details on selecting the value of $c$ and the estimation of the lower bound on $R$ in Section VI.
\end{remark}
\begin{remark}
In Theorem 1, the geometric convergence rate of DGM-BB-C is established given that $\frac{1}{n}\sum^{n}_{i=1}\alpha_{k}^{i} \leq\frac{2}{L}-\frac{\mu}{L^2}$ for all $k$, where $\frac{2}{L}-\frac{\mu}{L^2}\geq\frac{2}{\mu+L}$. It means that the average of the step sizes of all agents for all iterations has an upper bound which is no less than $\frac{2}{\mu+L}$. For every agent, the step size is automatically computed while running the algorithm, and it is no less than $\frac{1}{L}$ and does exceed the upper bound $\frac{1}{\mu}$. Hence, the step sizes of DGM-BB-C are larger than the scalar $\frac{1}{2L}$ or $\frac{1}{L}$ requested by these methods \cite{NedicDIGing17,Qu16,NedicATC17,Lu18,Xu18}. In Section VI, we demonstrate that the BB step sizes generated by our algorithm in numerical experiments meet their theoretical bounds.
\end{remark}
\section{Auxiliary Relations}
In this section, we provide several basic relations, which are prepared for the proof of Lemma 2 in our later analysis. An iterative equation that governs the average sequence  $\{\overline{y}_k\}$ to asymptotically track the average of local gradients is derived in Lemma 5. Lemma 6 presents some inequalities that are obtained directly from Assumption 1. Lemma 7 is a standard result in optimization theory, which states that the distance to optimizer shrinks by at least a fixed ratio in the centralized gradient method for a smooth and strongly convex function.
\begin{lemma}
There holds $\overline{y}_k=\frac{1}{n}\mathbf{1}^T\nabla \mathbf{f}(\mathbf{x}_k)$ for all $k\geq 0$.
\end{lemma}
 \begin{IEEEproof}
 Since $W$ is doubly stochastic satisfying $\mathbf{1}^TW=\mathbf{1}^T$, we obtain that
 \begin{IEEEeqnarray}{rCl}\label{41}
\overline{y}_k&=&\frac{1}{n}\mathbf{1}^TW^{R}[\mathbf{y}_{k-1}+\nabla \mathbf{f}(\mathbf{x}_{k})-\nabla \mathbf{f}(\mathbf{x}_{k-1})] \nonumber\\
&=&\overline{y}_{k-1}+\frac{1}{n}\mathbf{1}^T\nabla \mathbf{f}(\mathbf{x}_{k})-\frac{1}{n}\mathbf{1}^T\nabla \mathbf{f}(\mathbf{x}_{k-1}).\nonumber
\end{IEEEeqnarray}
Do this recursively, we have $\overline{y}_k=
\overline{y}_{0}+\frac{1}{n}\mathbf{1}^T\nabla \mathbf{f}(\mathbf{x}_{k})-\frac{1}{n}\mathbf{1}^T\nabla \mathbf{f}(\mathbf{x}_0).$ Since $\mathbf{y}_0=\nabla \mathbf{f}(\mathbf{x}_0)$, one can get $\overline{y}_0=\frac{1}{n}\mathbf{1}^T\nabla \mathbf{f}(\mathbf{x}_0)$. The proof is thus completed.
\end{IEEEproof}
\begin{lemma}
Under Assumption 1, for all $k\geq 0$, the following inequalities hold:\\
      $(i)$~~~$\|\nabla \mathbf{f}(\mathbf{x}_{k+1})-\nabla \mathbf{f}(\mathbf{x}_{k})\|\leq L\|\mathbf{x}_{k+1}-\mathbf{x}_{k}\|$ ;\\
      $(ii)$~~$\|\frac{1}{n}\mathbf{1}^T\nabla \mathbf{f}(\mathbf{x}_{k+1})-\frac{1}{n}\mathbf{1}^T\nabla \mathbf{f}(\mathbf{x}_{k})\|\leq L\frac{1}{\sqrt{n}}\|\mathbf{x}_{k+1}-\mathbf{x}_{k}\|$ ;\\
      $(iii)$~$\|\frac{1}{n}\mathbf{1}^T\nabla \mathbf{f}(\mathbf{x}_k)-\nabla f(\overline{x}_{k})\|\leq L\frac{1}{\sqrt{n}}\|\mathbf{x}_k-\mathbf{1}\overline{x}_{k}\|$.
      \end{lemma}
 \begin{IEEEproof}
      The results follow from the proof of Lemma 8 in \cite{Qu16}.
 \end{IEEEproof}
 \begin{lemma}
 For any $u\in \mathbb{R}^{1\times p}$, assume that function $h(u)$ is $\mu$-strongly convex and has Lipschitz continuous gradient with constant $L$. Let $u^*\in \mathbb{R}^{1\times p}$ be the global minimum of $h$ and $0<\alpha<\frac{2}{L}$. Then we have
\begin{eqnarray}\label{42}
     \|u-\alpha\nabla h(u)-u^*\| \leq \lambda_{h}\|u-u^*\|,
      \end{eqnarray}
     where $\lambda_{h}=\max\{|1-\alpha\mu|, |1-\alpha L|\}$.
     \end{lemma}
     \begin{IEEEproof}
     See Lemma 10 of \cite{Qu16} for reference.
     \end{IEEEproof}
\section{Convergence Analysis}
We now provide the proof of Lemma 2 in this section. We will bound $\|\mathbf{x}_{k+1}-\mathbf{1}\overline{x}_{k+1}\|, \|\mathbf{y}_{k+1}-\mathbf{1}\overline{y}_{k+1}\|$ and $\|\overline{x}_{k+1}-(x^*)^T\|$ in terms of linear combinations of their past values, in which way we establish a linear system of inequalities.
\begin{IEEEproof}
Step 1: Bound $\|\mathbf{x}_{k+1}-\mathbf{1}\overline{x}_{k+1}\|$.

From $(\ref{13})$, it follows that
\begin{IEEEeqnarray}{rCl}\label{43}
&&\|\mathbf{x}_{k+1}-\mathbf{1}\overline{x}_{k+1}\|\nonumber\\
&=&\|\mathbf{x}_{k+1}-\frac{1}{n}\mathbf{1}\mathbf{1}^T\mathbf{x}_{k+1}\| \nonumber\\
&=&\|(I-\frac{1}{n}\mathbf{1}\mathbf{1}^T)W^{R}[\mathbf{x}_{k}-D_{k}\mathbf{y}_{k}]\|\nonumber\\
&\leq&\|(I-\frac{1}{n}\mathbf{1}\mathbf{1}^T)W^{R}\mathbf{x}_{k}\|+\|(I-\frac{1}{n}\mathbf{1}\mathbf{1}^T)W^{R}D_{k}\mathbf{y}_{k}\|.
\end{IEEEeqnarray}
Using $(I-\frac{1}{n}\mathbf{1}\mathbf{1}^T)W=(W-\frac{1}{n}\mathbf{1}\mathbf{1}^T)(I-\frac{1}{n}\mathbf{1}\mathbf{1}^T)$ and by the definition of $\delta$, it follows from $(\ref{43})$ that
 \begin{IEEEeqnarray}{rCl}\label{44}
\|\mathbf{x}_{k+1}-\mathbf{1}\overline{x}_{k+1}\|&\leq& \|(W-\frac{1}{n}\mathbf{1}\mathbf{1}^T)^{R}(I-\frac{1}{n}\mathbf{1}\mathbf{1}^T)\mathbf{x}_{k}\|+\|(W-\frac{1}{n}\mathbf{1}\mathbf{1}^T)^{R}
D_{k}\mathbf{y}_{k}\| \nonumber\\
&\leq& \delta^{R}\|\mathbf{x}_{k}-\mathbf{1}\overline{x}_{k}\|+\delta^{R}\alpha_{\max}\|\mathbf{y}_{k}\|.
\end{IEEEeqnarray}
By Lemmas 5 and 6 and due to $\mathbf{1}^T\nabla \mathbf{f}(\mathbf{x}^*)=0$, we have
\begin{IEEEeqnarray}{rCl}\label{45}
\|\mathbf{y}_{k}\|&=&\|\mathbf{y}_{k}-\frac{1}{n}\mathbf{1}\mathbf{1}^T\mathbf{y}_{k}+\frac{1}{n}\mathbf{1}\mathbf{1}^T\mathbf{y}_{k}
-\mathbf{1}\nabla f(\overline{x}_{k})+\mathbf{1}\nabla f(\overline{x}_{k})-\frac{1}{n}\mathbf{1}\mathbf{1}^T\nabla \mathbf{f}(\mathbf{x}^*)\| \nonumber\\
&\leq&  \|\mathbf{y}_{k}-\mathbf{1}\overline{y}_{k}\|+\|\mathbf{1}[\frac{1}{n}\mathbf{1}^T\nabla \mathbf{f}(\mathbf{x}_{k})-\nabla f(\overline{x}_{k})]\|+\|\mathbf{1}[\nabla f(\overline{x}_{k})-\frac{1}{n}\mathbf{1}^T\nabla \mathbf{f}(\mathbf{x}^*)]\| \nonumber \\
&\leq& \|\mathbf{y}_{k}-\mathbf{1}\overline{y}_{k}\|+L\|\mathbf{x}_{k}-\mathbf{1}\overline{x}_{k}\|
+L\sqrt{n}\|\overline{x}_{k}-(x^*)^T\|.
\end{IEEEeqnarray}
Substituting $(\ref{45})$ into $(\ref{44})$ yields
\begin{IEEEeqnarray}{rCl}\label{46}
\|\mathbf{x}_{k+1}-\mathbf{1}\overline{x}_{k+1}\| \leq (\delta^{R}+\delta^{R}L\alpha_{\max})\|\mathbf{x}_{k}-\mathbf{1}\overline{x}_{k}\| \nonumber\\
+\delta^{R}\alpha_{\max}\|\mathbf{y}_{k}-\mathbf{1}\overline{y}_{k}\|
+\delta^{R}L\sqrt{n}\alpha_{\max}\|\overline{x}_{k}-(x^*)^T\|.
\end{IEEEeqnarray}

Step 2: Bound $\|\mathbf{y}_{k+1}-\mathbf{1}\overline{y}_{k+1}\|$.

Noticing the similarity between $(\ref{13})$ and $(\ref{14})$, it follows from $(\ref{14})$ that
\begin{IEEEeqnarray}{rCl}\label{47}
&&\|\mathbf{y}_{k+1}-\mathbf{1}\overline{y}_{k+1}\|\nonumber\\
&\leq& \delta^R\|\mathbf{y}_{k}-\mathbf{1}\overline{y}_{k}\|+\delta^R\|\nabla \mathbf{f}(\mathbf{x}_{k+1})-\nabla \mathbf{f}(\mathbf{x}_{k})\| \nonumber\\
&\leq& \delta^R\|\mathbf{y}_{k}-\mathbf{1}\overline{y}_{k}\|+\delta^RL\|\mathbf{x}_{k+1}-\mathbf{x}_{k}\|,
\end{IEEEeqnarray}
where we have used Lemma 6($i$) in the last inequality $(\ref{47})$.
Let us look into the last term in  $(\ref{47})$, we have
\begin{IEEEeqnarray}{rCl}\label{48}
\|\mathbf{x}_{k+1}-\mathbf{x}_{k}\|&=&\|W^{R}[\mathbf{x}_{k}-D_{k}\mathbf{y}_{k}]-\mathbf{x}_{k}\|\nonumber \\
&=&\|(W^{R}-I)[\mathbf{x}_{k}-\mathbf{1}\overline{x}_{k}]-W^{R}D_{k}\mathbf{y}_{k}\|\nonumber \\
&<&2\|\mathbf{x}_{k}-\mathbf{1}\overline{x}_{k}\|+\alpha_{\max}\|\mathbf{y}_{k}\|.
\end{IEEEeqnarray}
 By combining $(\ref{45})$, $(\ref{47})$ with $(\ref{48})$, one has
 \begin{IEEEeqnarray}{rCl}\label{49}
\|\mathbf{y}_{k+1}-\mathbf{1}\overline{y}_{k+1}\|&\leq& (2\delta^RL+\delta^RL^2\alpha_{\max})\|\mathbf{x}_{k}-\mathbf{1}\overline{x}_{k}\|\nonumber\\
&+&(\delta^R+\delta^RL\alpha_{\max})\|\mathbf{y}_{k}-\mathbf{1}\overline{y}_{k}\|
\nonumber\\
&+&\delta^RL^2\sqrt{n}\alpha_{\max}\|\overline{x}_{k}-(x^*)^T\|.
\end{IEEEeqnarray}

Step 3: Bound $\|\overline{x}_{k+1}-(x^*)^T\|$.

 Since $W$ is doubly stochastic, we have $\mathbf{1}^TW=\mathbf{1}^T$. Taking the average of $(\ref{13})$ over $i$ gives us
\begin{IEEEeqnarray}{rCl}\label{50}
\overline{x}_{k+1}&=&\frac{1}{n}\mathbf{1}^TW^{R}[\mathbf{x}_{k}-D_{k}\mathbf{y}_{k}]\nonumber\\
&=&\overline{x}_{k}-\frac{1}{n}\mathbf{1}^TD_{k}\mathbf{y}_{{k}}.
\end{IEEEeqnarray}
It follows from  $(\ref{50})$ that
\begin{IEEEeqnarray}{rCl}\label{51}
&&\|\overline{x}_{k+1}-(x^*)^T\|\nonumber\\
&=&\|\overline{x}_{k}-\frac{1}{n}\mathbf{1}^TD_{k}\mathbf{y}_{{k}}-(x^*)^T\|\nonumber\\
&=&\|\overline{x}_{k}-\frac{1}{n}\mathbf{1}^TD_{k}\mathbf{y}_{{k}}+\frac{1}{n}\mathbf{1}^T(D_{k}-D_{k})\mathbf{1}\overline{y}_{{k}}-(x^*)^T\| \nonumber\\
&\leq& \|\overline{x}_{k}-\frac{1}{n}\mathbf{1}^TD_{k}\mathbf{1}\overline{y}_{{k}}-(x^*)^T\|+\|\frac{1}{n}\mathbf{1}^TD_{k}(\mathbf{y}_{{k}}-\mathbf{1}\overline{y}_{{k}})\|\nonumber\\
&\leq& \|\overline{x}_{k}-\frac{1}{n}\mathbf{1}^TD_{k}\mathbf{1}\overline{y}_{{k}}-(x^*)^T\|+\frac{\alpha_{\max}}{\sqrt{n}}\|\mathbf{y}_{{k}}-\mathbf{1}\overline{y}_{{k}}\|.
\end{IEEEeqnarray}
Consider the first term in $(\ref{51})$, by Lemmas 5 and 6, we derive that
 \begin{IEEEeqnarray}{rCl}\label{52}
 &&\|\overline{x}_{k}-\frac{1}{n}\mathbf{1}^TD_{k}\mathbf{1}\overline{y}_{{k}}-(x^*)^T\|\nonumber\\
 &=&\|\overline{x}_{k}-(\frac{1}{n}\sum^{n}_{i=1}\alpha_{k}^{i})\overline{y}_{{k}}-(x^*)^T\|\nonumber\\
 &\leq& \|\overline{x}_{k}-(\frac{1}{n}\sum^{n}_{i=1}\alpha_{k}^{i})\nabla f(\overline{x}_{k})-(x^*)^T\|
 +\frac{1}{n}\sum^{n}_{i=1}\alpha_{k}^{i}\|\frac{1}{n}\mathbf{1}^T\nabla \mathbf{f}(\mathbf{x}_k)-\nabla f(\overline{x}_{k})\| \nonumber \\
 &\leq& \|\overline{x}_{k}-(\frac{1}{n}\sum^{n}_{i=1}\alpha_{k}^{i})\nabla f(\overline{x}_{k})-(x^*)^T\|
 +(\frac{1}{n}\sum^{n}_{i=1}\alpha_{k}^{i})\frac{L}{\sqrt{n}}\|\mathbf{x}_{k}-\mathbf{1}\overline{x}_{k}\|.
\end{IEEEeqnarray}
If $\frac{1}{n}\sum^{n}_{i=1}\alpha_{k}^{i}<\frac{2}{L}$, according to Lemma 7, we have
\begin{IEEEeqnarray}{rCl}\label{53}
 \|\overline{x}_{k}\!-\!(\frac{1}{n}\sum^{n}_{i=1}\alpha_{k}^{i})\nabla f(\overline{x}_{k})\!-\!(x^*)^T\|\leq \lambda\|\overline{x}_{k}\!-\!(x^*)^T\|,
 \end{IEEEeqnarray}
 where $\lambda=\max\{|1-\frac{\mu}{n}\sum^{n}_{i=1}\alpha_{k}^{i}|, |1-\frac{L}{n}\sum^{n}_{i=1}\alpha_{k}^{i}|\}$.\\
 Recalling the definition of $\alpha_{\max}$ and by Lemma 1, we get $\frac{1}{n}\sum^{n}_{i=1}\alpha_{k}^{i}\leq \alpha_{\max}$. Combining $(\ref{52})$ with $(\ref{53})$, we obtain
\begin{IEEEeqnarray}{rCl}\label{54}
\|\overline{x}_{k}-\frac{1}{n}\mathbf{1}^TD_{k}\mathbf{1}\overline{y}_{{k}}-(x^*)^T\|
 \leq \lambda\|\overline{x}_{k}-(x^*)^T\|
 +\frac{L\alpha_{\max}}{\sqrt{n}}\|\mathbf{x}_{k}-\mathbf{1}\overline{x}_{k}\|.
\end{IEEEeqnarray}
By substituting $(\ref{54})$ into $(\ref{51})$, one gets
\begin{IEEEeqnarray}{rCl}\label{55}
\|\overline{x}_{k+1}-(x^*)^T\|\leq \frac{L\alpha_{\max}}{\sqrt{n}}\|\mathbf{x}_{k}-\mathbf{1}\overline{x}_{k}\|\nonumber\\
+\frac{\alpha_{\max}}{\sqrt{n}}\|\mathbf{y}_{{k}}-\mathbf{1}\overline{y}_{{k}}\|+ \lambda\|\overline{x}_{k}-(x^*)^T\|.
\end{IEEEeqnarray}
The proof is completed.
\end{IEEEproof}
\section{Numerical Experiments}
In this section, we analyze the performance of DGM-BB-C and illustrate our theoretical findings. Our numerical experiments are based on a distributed least squares problem over undirect graph which is generated by using the Erd\H{o}s-R\'{e}nyi model with connectivity ratio $r_c$ \cite{Erdos59}. We use the Metropolis constant edge weight matrix $W$ \cite{BoydW04}.

We consider a distributed sensing problem for solving an unknown signal $x\in\mathbb{R}^p$ \cite{Li17}. Each agent $i\in \{1,2,...,n\} $ holds its own measurement equation, $y_i=M_ix+e_i$, where $y_i\in \mathbb{R}^{m_i}$ and $M_i\in \mathbb{R}^{m_i\times p}$ are measured data, $e_i\in \mathbb{R}^{m_i}$ is unknown noise. We apply the least squares loss and try to solve
\begin{IEEEeqnarray}{rCl}\label{57}
        \min\limits_{x\in \mathbb{R}^p} f(x)=\frac{1}{n}\sum^{n}_{i=1}\frac{1}{2}\|M_{i}x-y_{i}\|^2_2.
\end{IEEEeqnarray}

In our experiment, we set $n=200, m_i=20, p=10$. In order to ensure that the objective function $f(x)$ is strongly convex, we choose $M_i$ such that the Lipschitz constant of $\nabla f(x)$ satisfies $L=1$ and the strongly convex constant $\mu=0.5$. We initialize $x_0^i=0$ for every $i\in\mathcal{V}$ and choose the BB step size $\alpha^i_k$ computed by $(\ref{9})$ (the practical performance of $\alpha^i_k$ computed by $(\ref{8})$ is similar). We test the performance with $r_c=0.1$ and $r_c=0.3$, respectively. As pointed out by Remark 3, in order to ensure that the range of the largest step size given in $(\ref{16})$ is contained in the range given in $(\ref{30})$, we require that $\widehat{\alpha}>\frac{1}{\mu}$, which depends on $R$ and the adjustable parameter vector $c$. As shown in Theorem 1, DGM-BB-C converges geometrically to the optimal solution given that $\frac{1}{n}\sum^{n}_{i=1}\alpha_{k}^{i}\leq\frac{2}{L}-\frac{\mu}{L^2}$ for all $k$. To make the condition hold, it is sufficient to require $\overline{\alpha}_{\max}\leq\frac{2}{L}-\frac{\mu}{
L^2}$. Then, we can estimate theoretical bounds of $R$, $\alpha_{\max}$ and $\overline{\alpha}_{\max}$ and the value of $\widehat{\alpha}$ as follows:

$(i)$~~For $r_c=0.1$, $R\geq4$ with $c\!=\![0.9240,0.9889,0.6453
]^T$, $\alpha_{\max}\in [1,2]$, $\overline{\alpha}_{\max} \in[1,1.5]$ and $\widehat{\alpha}= 2.3853$.

$(ii)$~For $r_c=0.3$, $R\geq3$ with $c\!=\![ 0.9978, 1.0737, 0.6677]^T$, $\alpha_{\max}\in [1,2]$, $\overline{\alpha}_{\max} \in[1,1.5]$ and $\widehat{\alpha}= 2.2793$.

 We can easily observe that the number of inner consensus iterations $R$ can be moderate to guarantee that the range of the largest step size given in $(\ref{16})$ is contained in the range given in $(\ref{30})$, which is in agreement with the theory and Remark 3.

 Fig. 1 shows the performance of DGM-BB-C with different initial step size $\alpha^i_0$ for each $i$ when $r_c=0.1$. We can see that DGM-BB-C is not sensitive to the choice of $\alpha^i_0$. For different connectivity of network, there are similar results. Therefore, DGM-BB-C has very promising potential in practice because it generates the step sizes automatically while running the algorithm.
 \begin{figure}[!t]
\centering
\includegraphics[width=2.7in]{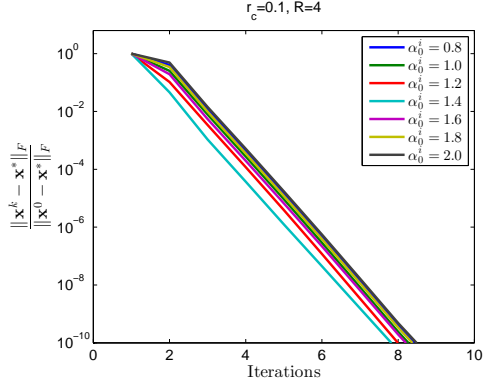}
\caption{The performance of DGM-BB-C with different $\alpha^i_0$.}
\label{fig:fig:1}
\end{figure}
For different connectivity of network, Fig. 2 shows the performance of DGM-BB-C with different number of inner consensus iterations $R$.  From Figs. 2$(a)$ and 2$(b)$, we can notice that DGM-BB-C with $R>1$ is faster than it with $R=1$, which illustrates the need for an additional consensus iterations. However, this does not mean that our algorithm performs better when $R$ is larger. It can be observed that our algorithm with $R=4$ performs best for less well-connected network, instead of $R=3$ for well-connected network.
\begin{figure}\label{fig:2}
\centering
  \subfloat[$r_c=0.1$]{
    \label{fig:2a} 
    \includegraphics[width=2.7in,height=2.3in]{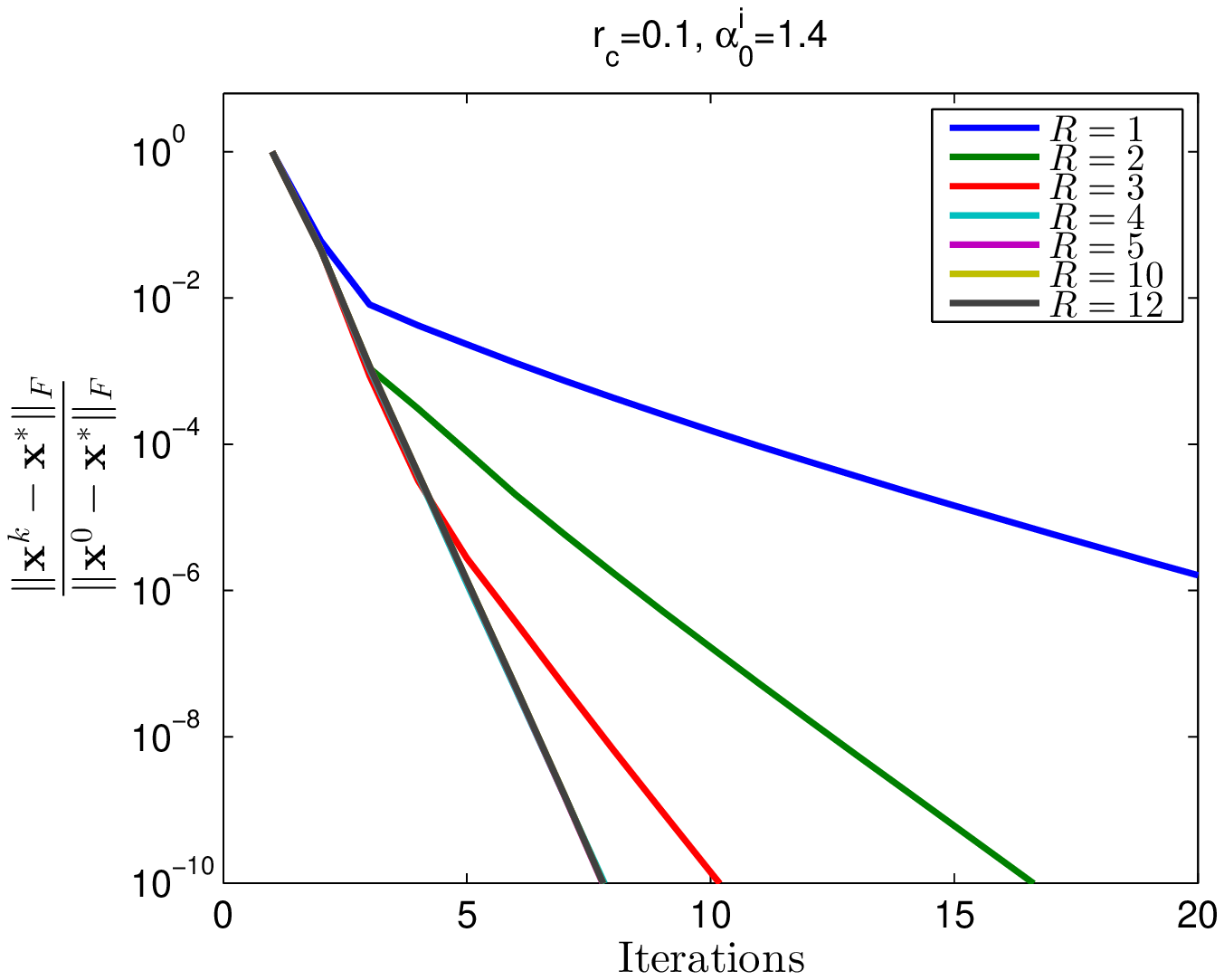}}
  \subfloat[$r_c=0.3$]{
    \label{fig:2b} 
    \includegraphics[width=2.7in,height=2.3in]{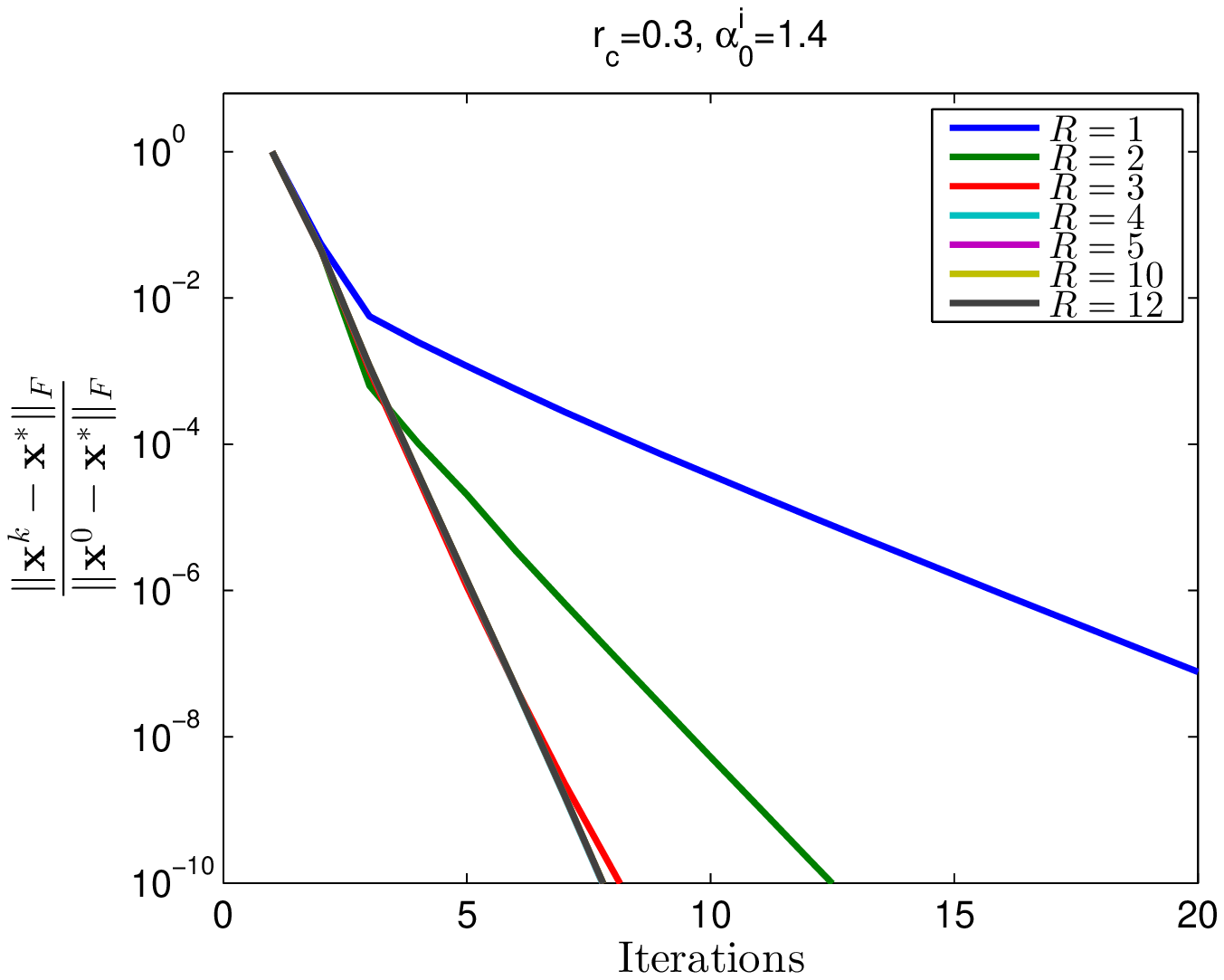}}
  \caption{The performance of DGM-BB-C with different $R$ on random network with $r_c=0.1$ and $r_c=0.3$.}
\end{figure}

In order to show that the proposed algorithm is more effective and comparable to some existing advanced methods, we do the corresponding numerical experiments. We compare the convergence rate of DGM-BB-C with several distributed algorithms, including DGD \cite{Nedic09}, NEAR-DGD$^{+}$ \cite{Berahas18}, EXTRA \cite{Shi14}, DIGing \cite{NedicDIGing17},  ATC-DIGing \cite{NedicATC17}. For ATC-DIGing, we firstly tune an optimized identical step size by hand for all agents, which is $\frac{1}{L}$ in our experiment and then perturb it by random variables satisfying the uniform distribution over interval $(0.6,1.2)$. For DGM-BB-C, we set $\alpha^i_0=1.4$ and $R=4$ for $r_c=0.1$ and $R=3$ for $r_c=0.3$. The step sizes employed in other algorithms are hand-optimized. We also analyze the convergence rate of DGM-C, which is a practical variant of DGM-BB-C. DGM-C uses the same step size rule as that of ATC-DIGing instead of using the BB step sizes.
\begin{figure}\label{fig:3}
\centering
  \subfloat[]{
    \label{fig:3a} 
    \includegraphics[width=2.7in,height=2.3in]{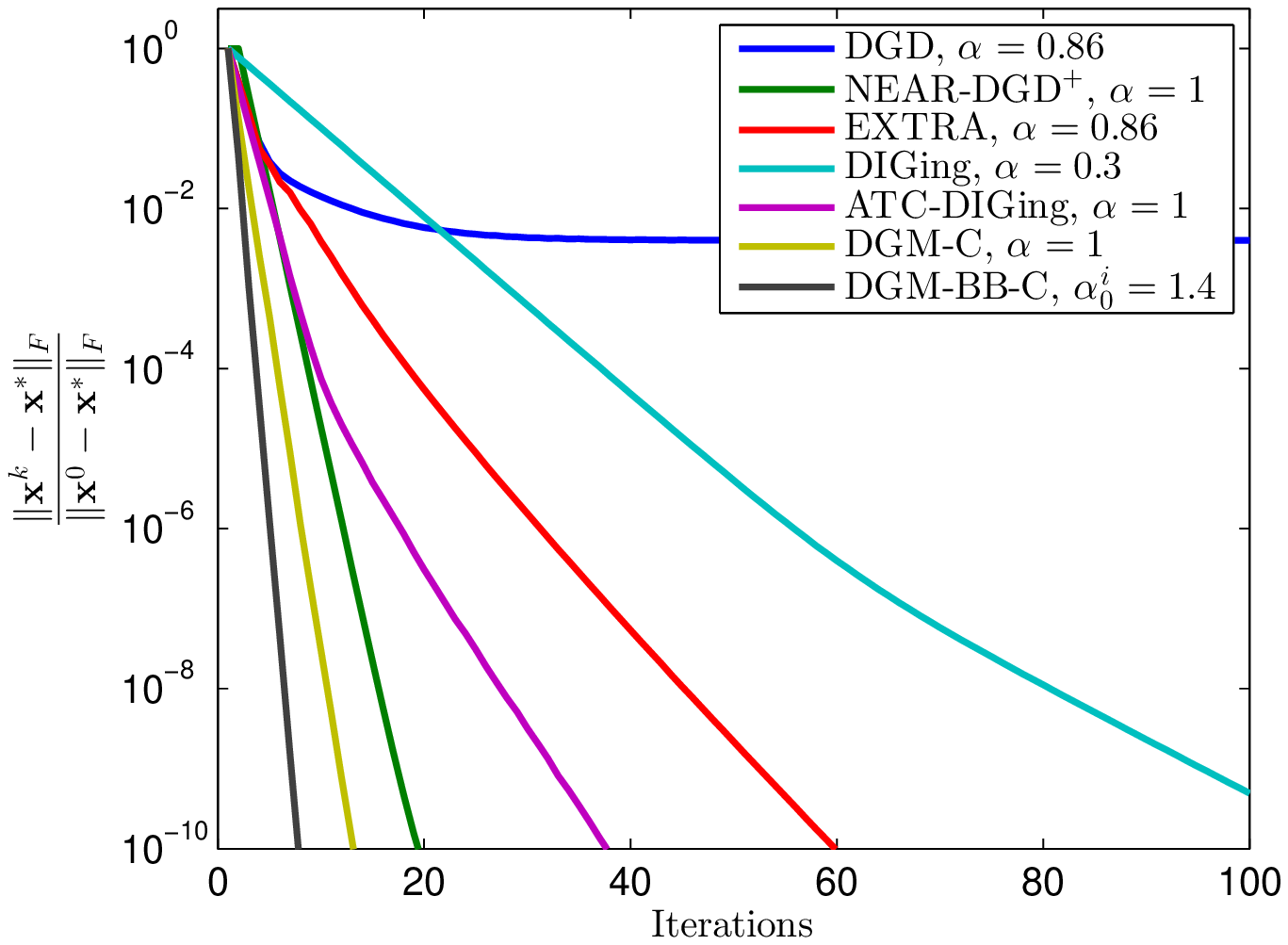}}
  \subfloat[]{
    \label{fig:3b} 
    \includegraphics[width=2.7in,height=2.3in]{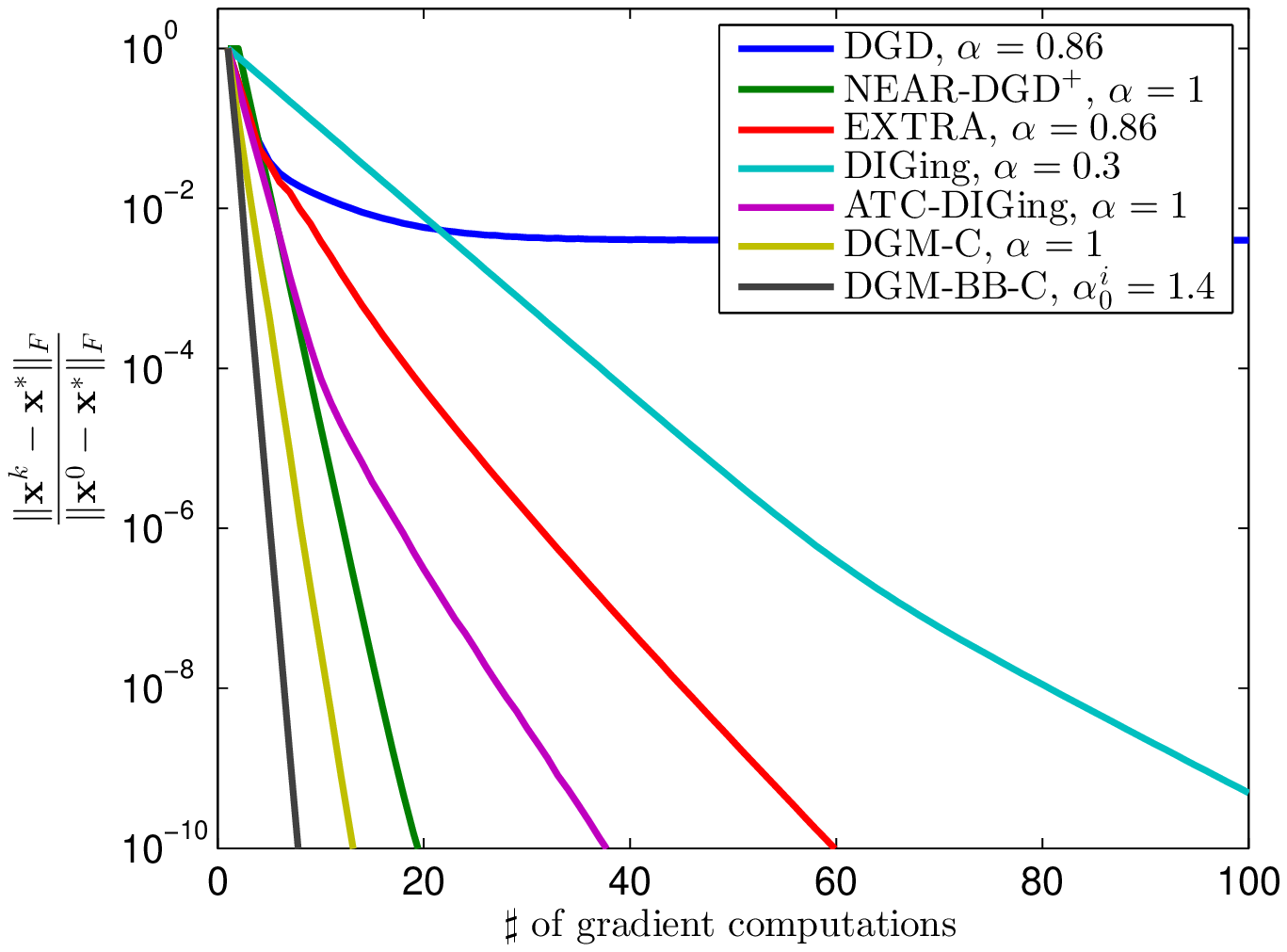}}\\
  \subfloat[]{
    \label{fig:3c} 
    \includegraphics[width=2.7in,height=2.3in]{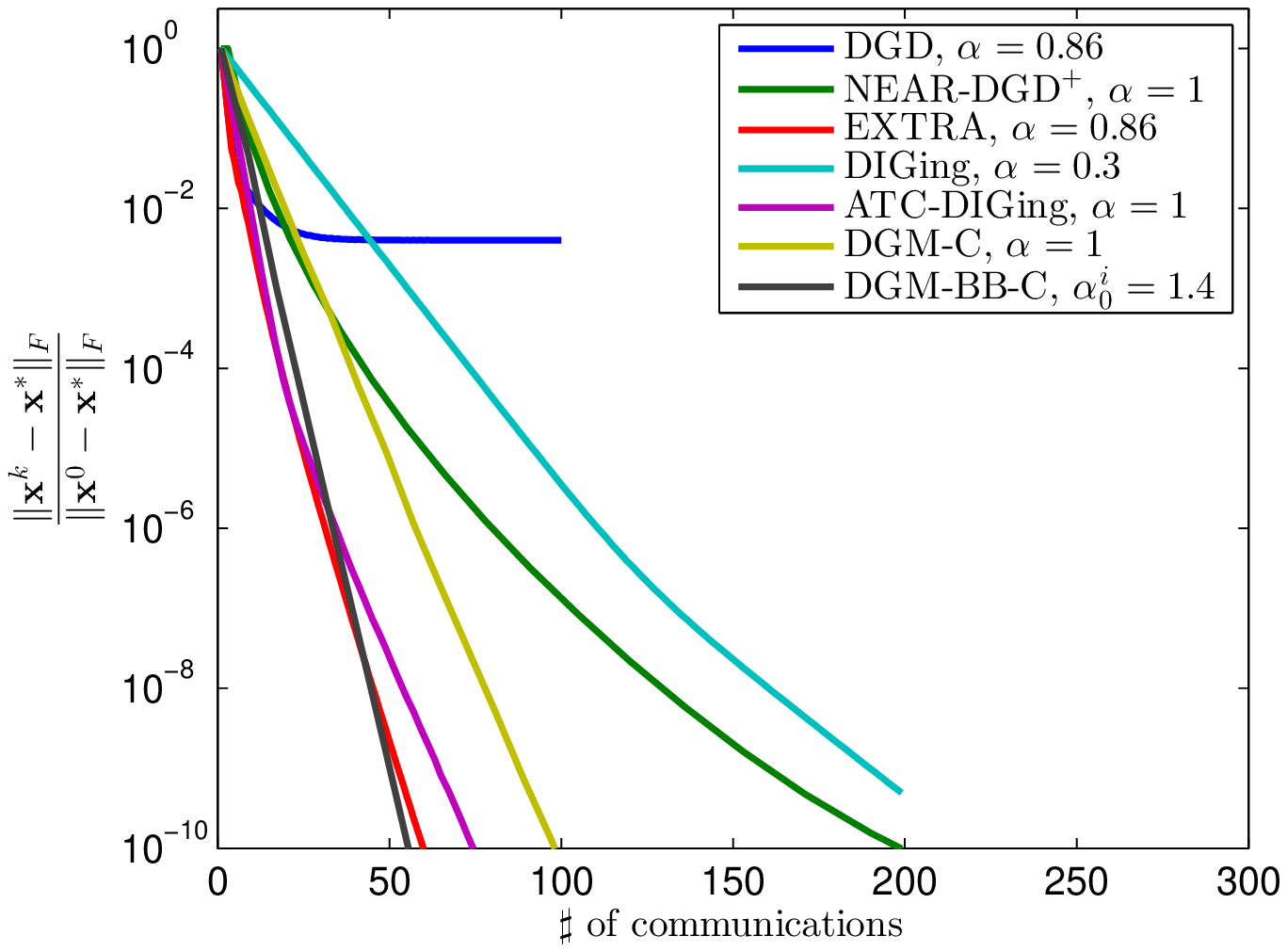}}
  \subfloat[]{
    \label{fig:3d} 
    \includegraphics[width=2.7in,height=2.3in]{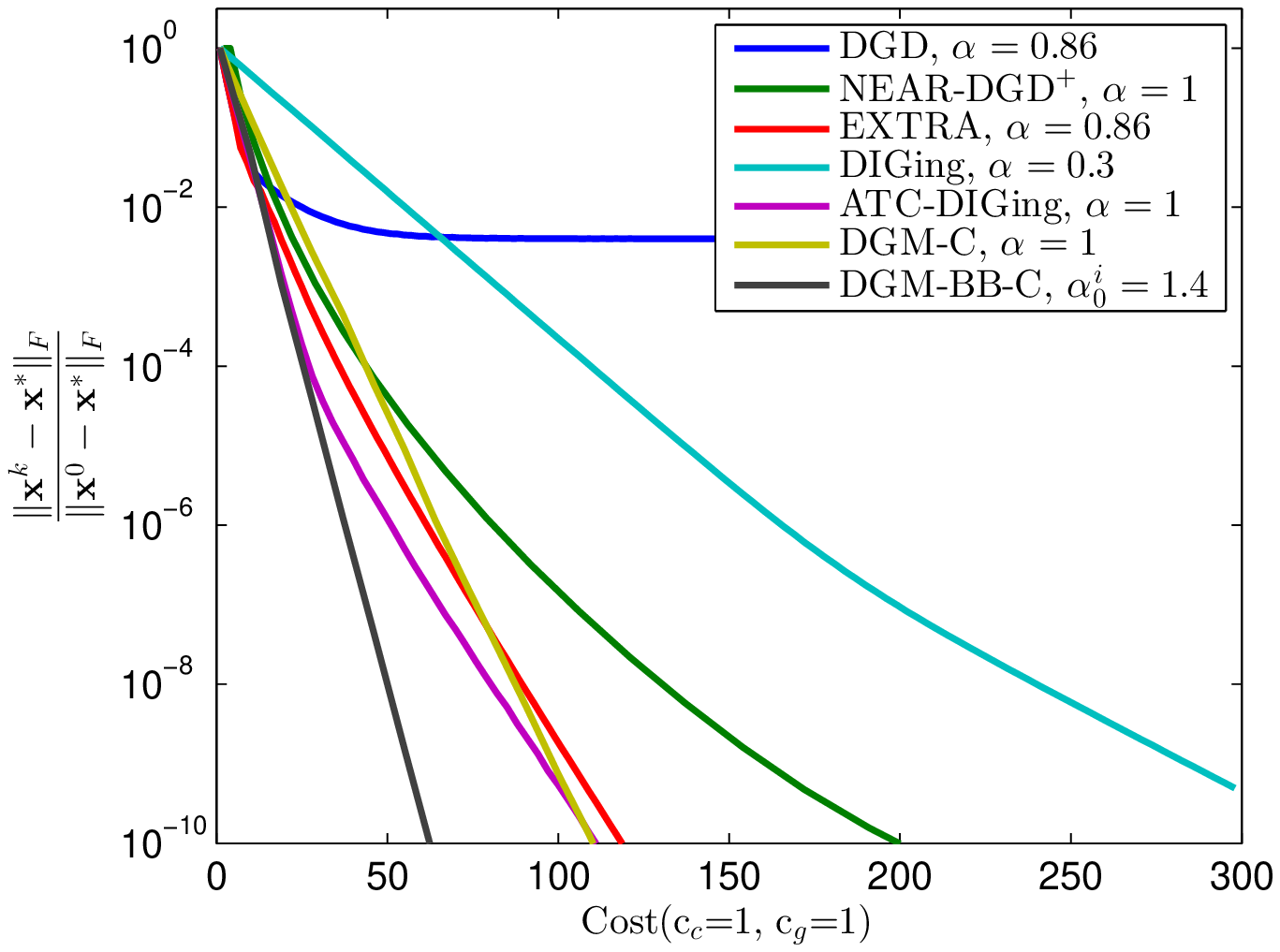}}
  \caption{Convergence rates comparison of different distributed algorithms with $r_c=0.1$.}
\end{figure}
\begin{figure}\label{fig:4}
\centering
  \subfloat[]{
    \label{fig:4a} 
    \includegraphics[width=2.7in,height=2.3in]{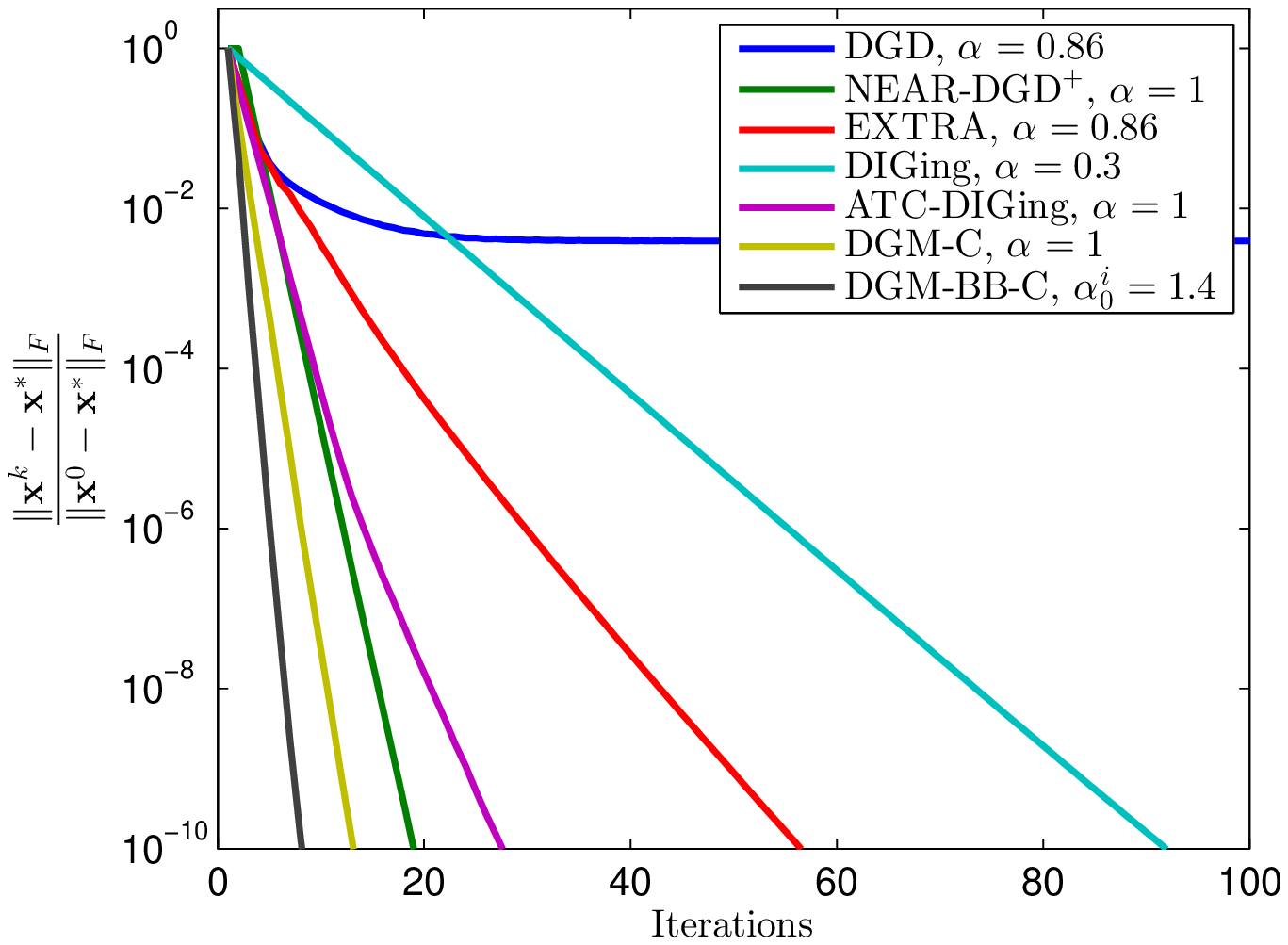}}
  \subfloat[]{
    \label{fig:4b} 
    \includegraphics[width=2.7in,height=2.3in]{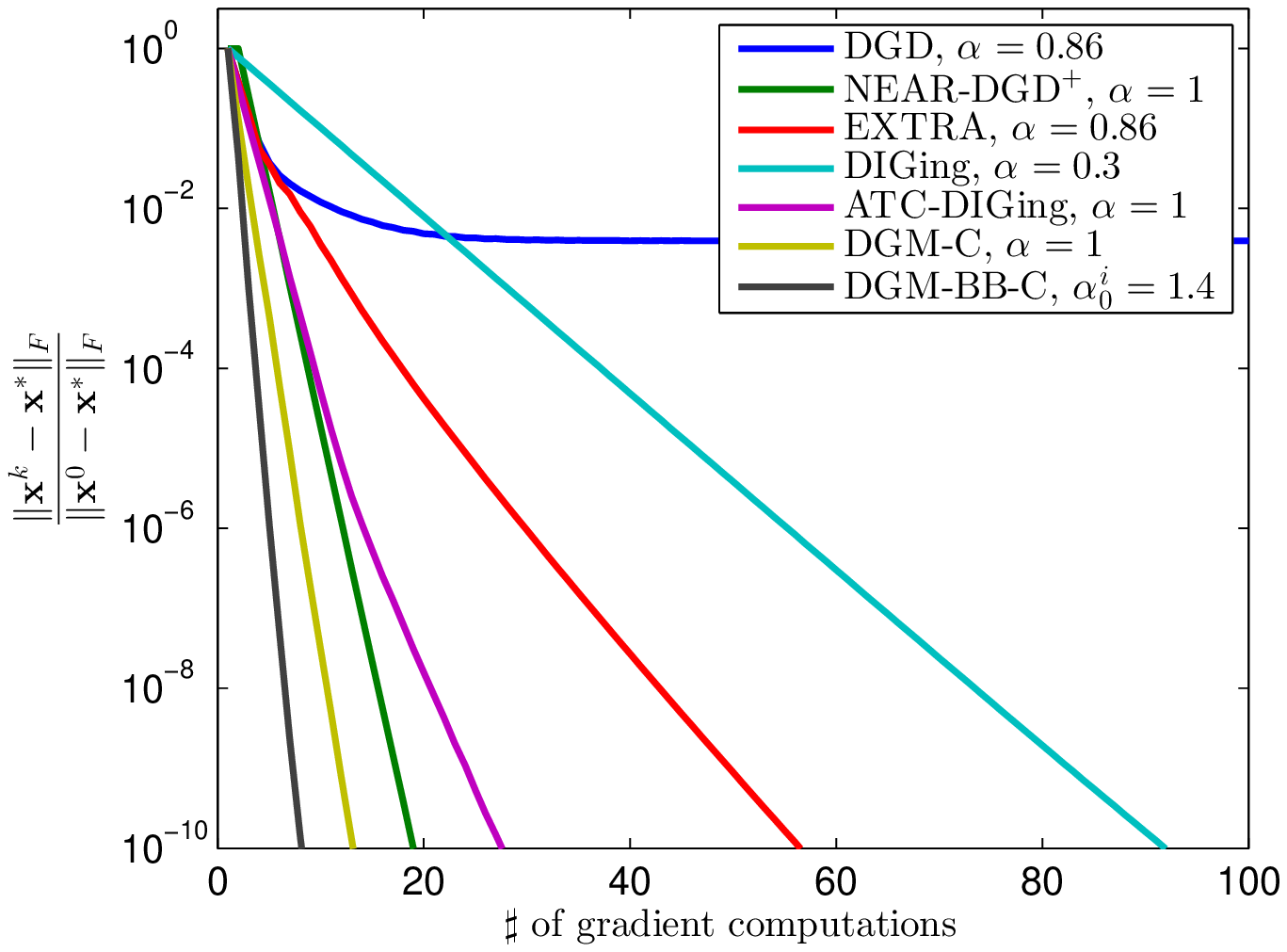}}\\
  \subfloat[]{
    \label{fig:4c} 
    \includegraphics[width=2.7in,height=2.3in]{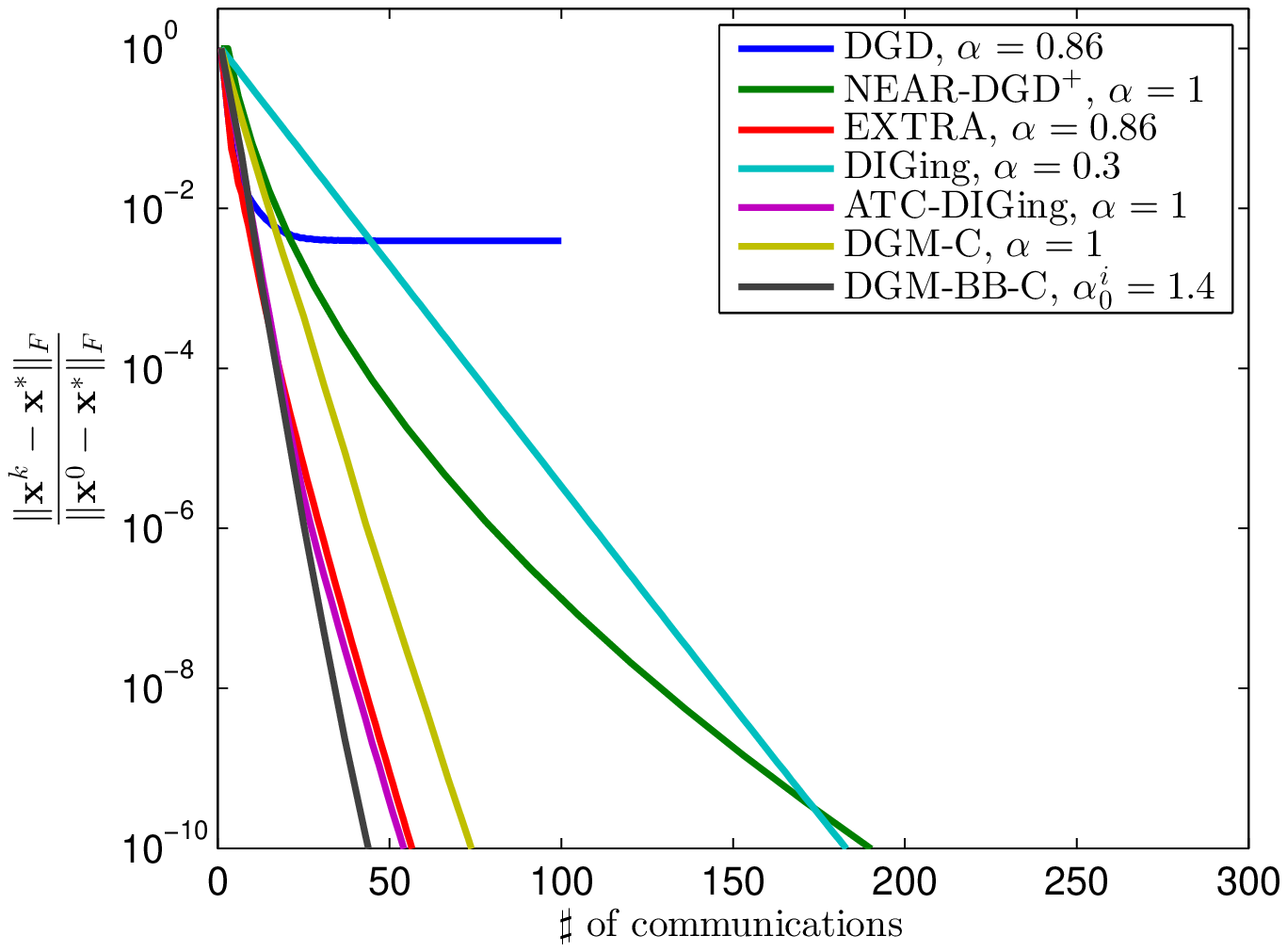}}
  \subfloat[]{
    \label{fig:4d} 
    \includegraphics[width=2.7in,height=2.3in]{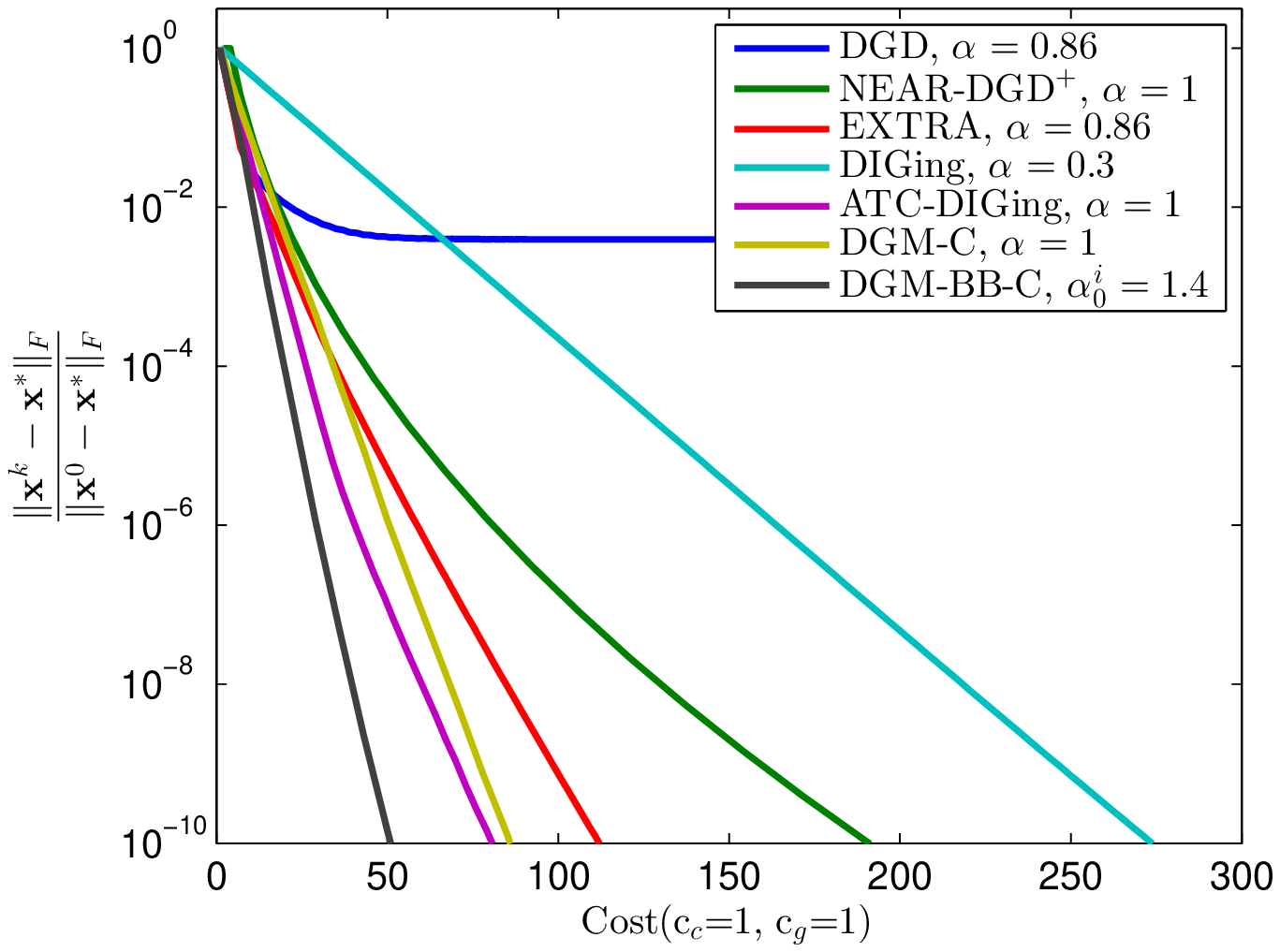}}
  \caption{Convergence rates comparison of different distributed algorithms with $r_c=0.3$.}
\end{figure}

In Figs. 3 and 4, we plot the relative error $\|\mathbf{x}_k-\mathbf{x}^*\| /\|\mathbf{x}_0-\mathbf{x}^*\|$ with respect to: $(i)$~ iterations, $(ii)$~ number of gradient computations, $(iii)$~number of communications, $(iv)$~cost (as described in \cite{Berahas18} with $c_c=c_g=1$). DGM-C performs better than ATC-DIGing which suggests that doing inner loops of consensus iterations can improve the performance of the algorithm. It follows from Figs. 3 and 4 that DGM-BB-C has the smallest number of iterations, gradient computation, and communication, and the lowest adaptive cost to reach an $\epsilon$-optimal solution because it allows to use larger step sizes and does not need to increase the number of consensus steps in practice and theory. Even though, at each iteration, having more communications, our algorithm has the lower computation cost and communication cost to reach an $\epsilon$-optimal solution, which means that DGM-BB-C achieves both the optimal computation and communication cost for distributed optimization. NEAR-DGD$^{+}$ seeks to the exact solution by increasing the number of consensus steps linearly with the iteration number, whereas DGM-BB-C can converge to the exact solution when the number of consensus steps stays constant.

 To further illustrate that our theoretical results are consistent with the numerical experimental results, we compute $\alpha_{\max}$ and $\overline{\alpha}_{\max}$ of the BB step sizes generated by our algorithm when the relative error reaches $10^{-10}$-accuracy. We also compute $\rho(G^{\alpha})$ because the convergence rate of DGM-BB-C is related to $\rho(G^{\alpha})$. These values are computed as follows:

 $(i)$~~For $r_c=0.1$, we have $\alpha_{\max}=1.8360$, $\overline{\alpha}_{\max}=1.3262$ and $\rho(G^{\alpha})=0.8713$.

 $(ii)$~For $r_c=0.3$, we have $\alpha_{\max}=1.8358$,  $\overline{\alpha}_{\max}=1.3308$ and $\rho(G^{\alpha})=0.8712$.

  We can see that $\alpha_{\max}$ and $\overline{\alpha}_{\max}$ are in the range of our theoretical value estimated above. We provide an estimate for the convergence rate.
\section{Conclusion}
In this paper, based on an adapt-then-combine variation of the dynamic average consensus approach and using multi-consensus inner loops, we propose a DGM-BB-C method. In contrast to the exiting distributed gradient methods, our method computes the step size for each agent automatically which only depends on its local information and is independent of the step sizes for other agents, and always admits the larger step sizes. Most importantly, for smooth and strongly convex objective functions, we have proved that DGM-BB-C converges geometrically to the optimal solution. We conduct numerical experiments on the distributed least squares problem, which has showed that our DGM-BB-C achieves both the optimal computation and communication cost for distributed optimization. DGM-BB-C can seek the exact solution both theoretically and empirically when the number of consensus steps stays constant. A possible topic in the future work is to extend our result to stochastic networks where we have to deal with asynchronous communication.


\ifCLASSOPTIONcaptionsoff
  \newpage
\fi


\begin{thebibliography}{10}
\providecommand{\url}[1]{#1}
\csname url@samestyle\endcsname
\providecommand{\newblock}{\relax}
\providecommand{\bibinfo}[2]{#2}
\providecommand{\BIBentrySTDinterwordspacing}{\spaceskip=0pt\relax}
\providecommand{\BIBentryALTinterwordstretchfactor}{4}
\providecommand{\BIBentryALTinterwordspacing}{\spaceskip=\fontdimen2\font plus
\BIBentryALTinterwordstretchfactor\fontdimen3\font minus
  \fontdimen4\font\relax}
\providecommand{\BIBforeignlanguage}[2]{{%
\expandafter\ifx\csname l@#1\endcsname\relax
\typeout{** WARNING: IEEEtran.bst: No hyphenation pattern has been}%
\typeout{** loaded for the language `#1'. Using the pattern for}%
\typeout{** the default language instead.}%
\else
\language=\csname l@#1\endcsname
\fi
#2}}
\providecommand{\BIBdecl}{\relax}
\BIBdecl

\bibitem{Cevher14}
V.~{Cevher}, S.~{Becker}, and M.~{Schmidt}, ``Convex optimization for big data:
  Scalable, randomized, and parallel algorithms for big data analytics,''
  \emph{IEEE Signal Process. Mag.}, vol.~31, no.~5, pp. 32--43, Sep. 2014.

\bibitem{Boyd11}
S.~{Boyd}, N.~{Parikh}, E.~{Chu}, B.~{Peleato}, and J.~{Eckstein},
  ``Distributed optimization and statistical learning via the alternating
  direction method of multipliers,'' \emph{Found. Trends Mach. Learn.}, vol.~3,
  no.~1, pp. 11--22, Jan. 2011.

\bibitem{Necoara08}
I.~{Necoara} and J.~A.~K. {Suykens}, ``Application of a smoothing technique to
  decomposition in convex optimization,'' \emph{IEEE Trans. Autom. Control},
  vol.~53, no.~11, pp. 2674--2679, Dec. 2008.

\bibitem{Bazerque10}
J.~A. {Bazerque} and G.~B. {Giannakis}, ``Distributed spectrum sensing for
  cognitive radio networks by exploiting sparsity,'' \emph{IEEE Trans. Signal
  Process.}, vol.~58, no.~3, pp. 1847--1862, Mar. 2010.

\bibitem{Olshevsky10}
A.~{Olshevsky}, ``Efficient information aggregation strategies for distributed
  control and signal processing,'' PhD thesis, Massachusetts Inst. Tech., 2010.

\bibitem{Ren06}
W.~{Ren}, ``Consensus based formation control strategies for multi-vehicle
  systems,'' in \emph{2006 Amer. Control Conf.}, Jun. 2006, pp. 4237--4242.

\bibitem{Pu16}
S.~{Pu}, A.~{Garcia}, and Z.~{Lin}, ``Noise reduction by swarming in social
  foraging,'' \emph{IEEE Trans. Autom. Control}, vol.~61, no.~12, pp.
  4007--4013, Dec. 2016.

\bibitem{Cohen17}
K.~{Cohen}, A.~{Nedi\'{c}}, and R.~{Srikant}, ``Distributed learning algorithms
  for spectrum sharing in spatial random access wireless networks,'' \emph{IEEE
  Trans. Autom. Control}, vol.~62, no.~6, pp. 2854--2869, Jun. 2017.

\bibitem{Ling10}
Q.~{Ling} and Z.~{Tian}, ``Decentralized sparse signal recovery for compressive
  sleeping wireless sensor networks,'' \emph{IEEE Trans. Signal Process.},
  vol.~58, no.~7, pp. 3816--3827, Jul. 2010.

\bibitem{Gan13}
L.~{Gan}, U.~{Topcu}, and S.~H. {Low}, ``Optimal decentralized protocol for
  electric vehicle charging,'' \emph{IEEE Trans. Power Syst.}, vol.~28, no.~2,
  pp. 940--951, May 2013.

\bibitem{Rabbat04}
M.~{Rabbat} and R.~{Nowak}, ``Distributed optimization in sensor networks,'' in
  \emph{Proc. 3rd int. symp. on Inform. process. sensor netw.}, Apr. 2004, pp.
  20--27.

\bibitem{Bertsekas83}
D.~{Bertsekas}, ``Distributed asynchronous computation of fixed points,''
  \emph{Math. Program.}, vol.~27, no.~1, pp. 107--120, 1983.

\bibitem{Tsitsiklis86}
J.~{Tsitsiklis}, D.~{Bertsekas}, and M.~{Athans}, ``Distributed asynchronous
  deterministic and stochastic gradient optimization algorithms,'' \emph{IEEE
  Trans. Autom. Control}, vol.~31, no.~9, pp. 803--812, Sep. 1986.

\bibitem{Nedic09}
A.~{Nedi\'{c}} and A.~{Ozdaglar}, ``Distributed subgradient methods for
  multi-agent optimization,'' \emph{IEEE Trans. Autom. Control}, vol.~54,
  no.~1, pp. 48--61, Jan. 2009.

\bibitem{JakoveticDNG14}
D.~{Jakoveti\'{c}}, J.~{Xavier}, and J.~M.~F. {Moura}, ``Fast distributed
  gradient methods,'' \emph{IEEE Trans. Autom. Control}, vol.~59, no.~5, pp.
  1131--1146, May 2014.

\bibitem{Yuan16}
K.~{Yuan}, Q.~{Ling}, and W.~{Yin}, ``On the convergence of decentralized
  gradient descent,'' \emph{SIAM J. Optim.}, vol.~26, no.~3, pp. 1835--1854,
  2016.

\bibitem{Terelius11}
H.~{Terelius}, U.~{Topcu}, and R.~{Murray}, ``Decentralized multi-agent
  optimization via dual decomposition,'' in \emph{18th IFAC Proceedings
  Volumes}, vol.~44, no.~1, Aug. 2011, pp. 11\,245--11\,251.

\bibitem{MotaADMM13}
J.~F.~C. {Mota}, J.~M.~F. {Xavier}, P.~M.~Q. {Aguiar}, and M.~{P\"{u}schel},
  ``D-admm: A communication-efficient distributed algorithm for separable
  optimization,'' \emph{IEEE Trans. Signal Process.}, vol.~61, no.~10, pp.
  2718--2723, May 2013.

\bibitem{ShiDADMM14}
W.~{Shi}, Q.~{Ling}, K.~{Yuan}, G.~{Wu}, and W.~{Yin}, ``On the linear
  convergence of the admm in decentralized consensus optimization,'' \emph{IEEE
  Trans. Signal Process.}, vol.~62, no.~7, pp. 1750--1761, Apr. 2014.

\bibitem{Mokhtari16}
A.~{Mokhtari}, W.~{Shi}, Q.~{Ling}, and A.~{Ribeiro}, ``A decentralized
  second-order method with exact linear convergence rate for consensus
  optimization,'' \emph{IEEE Trans. Signal Inform. Process. Netw.}, vol.~2,
  no.~4, pp. 507--522, Dec. 2016.

\bibitem{Eisen17}
M.~{Eisen}, A.~{Mokhtari}, and A.~{Ribeiro}, ``Decentralized quasi-newton
  methods,'' \emph{IEEE Trans. Signal Process.}, vol.~65, no.~10, pp.
  2613--2628, May 2017.

\bibitem{Berahas18}
A.~{Berahas}, R.~{Bollapragada}, N.~S. {Keskar}, and E.~{Wei}, ``Balancing
  communication and computation in distributed optimization,'' \emph{IEEE
  Trans. Autom. Control}, pp. 1--1, 2018.

\bibitem{Shi14}
W.~Shi, Q.~Ling, G.~Wu, and W.~Yin, ``Extra: An exact first-order algorithm for
  decentralized consensus optimization,'' \emph{SIAM J. Optim.}, vol.~25,
  no.~2, pp. 944--966, 2015.

\bibitem{XiDEXTRA17}
C.~{Xi} and U.~A. {Khan}, ``Dextra: A fast algorithm for optimization over
  directed graphs,'' \emph{IEEE Trans. Autom. Control}, vol.~62, no.~10, pp.
  4980--4993, Oct. 2017.

\bibitem{XuPhD16}
J.~{Xu}, ``Augmented distributed optimization for networked systems,'' PhD
  thesis, Nanyang Tech. Univ., 2016.

\bibitem{Xu15}
J.~{Xu}, S.~{Zhu}, Y.~C. {Soh}, and L.~{Xie}, ``Augmented distributed gradient
  methods for multi-agent optimization under uncoordinated constant
  stepsizes,'' in \emph{Proc. 54th IEEE Conf. Decis. Control (CDC)}, Dec. 2015,
  pp. 2055--2060.

\bibitem{NedicDIGing17}
A.~{Nedic}, A.~{Olshevsky}, and W.~{Shi}, ``Achieving geometric convergence for
  distributed optimization over time-varying graphs,'' \emph{SIAM J. Optim.},
  vol.~27, no.~4, pp. 2597--2633, 2017.

\bibitem{Qu16}
G.~{Qu} and N.~{Li}, ``Harnessing smoothness to accelerate distributed
  optimization,'' \emph{IEEE Trans. Control Netw. Syst.}, vol.~5, no.~3, pp.
  1245--1260, Sep. 2018.

\bibitem{NedicATC17}
A.~{Nedic}, A.~{Olshevsky}, W.~{Shi}, and C.~A. {Uribe}, ``Geometrically
  convergent distributed optimization with uncoordinated step-sizes,'' in
  \emph{2017 Amer. Control Conf. (ACC)}, May 2017, pp. 3950--3955.

\bibitem{Lu18}
Q.~{L{\"u}}, H.~{Li}, and D.~{Xia}, ``Geometrical convergence rate for
  distributed optimization with time-varying directed graphs and uncoordinated
  step-sizes,'' \emph{Inform. Sci.}, vol. 422, pp. 516--530, 2018.

\bibitem{Xu18}
J.~{Xu}, S.~{Zhu}, Y.~C. {Soh}, and L.~{Xie}, ``Convergence of asynchronous
  distributed gradient methods over stochastic networks,'' \emph{IEEE Trans.
  Autom. Control}, vol.~63, no.~2, pp. 434--448, Feb. 2018.

\bibitem{Xin19}
R.~{Xin}, C.~{Xi}, and U.~A. {Khan}, ``Frost--fast row-stochastic optimization
  with uncoordinated step-sizes,'' \emph{J. Adv. Signal Process.}, vol. 2019,
  no.~1, 2019.

\bibitem{Zhu10}
M.~{Zhu} and S.~{Martinez}, ``Discrete-time dynamic average consensus,''
  \emph{Automatica}, vol.~46, no.~2, pp. 322--329, 2010.

\bibitem{BB88}
J.~{Barzilai} and J.~M. {Borwein}, ``Two-point step size gradient methods,''
  \emph{IMA J. Numer. Anal.}, vol.~8, no.~1, pp. 141--148, 1988.

\bibitem{Raydan97}
M.~{Raydan}, ``The barzilai and borwein gradient method for the large scale
  unconstrained minimization problem,'' \emph{SIAM J. Optim.}, vol.~7, no.~1,
  pp. 26--33, 1997.

\bibitem{Dai05}
Y.~H. {Dai} and R.~{Fletcher}, ``Projected barzilai-borwein methods for
  large-scale box-constrained quadratic programming,'' \emph{Numer. Math.},
  vol. 100, no.~1, pp. 21--47, 2005.

\bibitem{Dai06}
Y.~H. {Dai}, W.~W. {Hager}, K.~{Schittkowski}, and H.~C. {Zhang}, ``The cyclic
  barzilai-borwein method for unconstrained optimization,'' \emph{IMA J. Numer.
  Anal.}, vol.~26, no.~3, pp. 604--627, 2006.

\bibitem{Dai13}
Y.~H. {Dai}, ``A new analysis on the barzilai-borwein gradient method,''
  \emph{J. oper. Res. Soc. China}, vol.~1, no.~2, pp. 187--198, 2013.

\bibitem{Dai19}
Y.~H. {Dai}, Y.~K. {Huang}, and X.~W. {Liu}, ``A family of spectral gradient
  methods for optimization,'' \emph{Comput. Optim. Appl.}, pp. 1--23, 2018.

\bibitem{Wang07}
Y.~F. {Wang} and S.~Q. {Ma}, ``Projected barzilai-borwein method for
  large-scale nonnegative image restoration,'' \emph{Inverse Probl. Sci. Eng.},
  vol.~15, no.~6, pp. 559--583, 2007.

\bibitem{Wright09}
S.~J. {Wright}, R.~D. {Nowak}, and M.~A.~T. {Figueiredo}, ``Sparse
  reconstruction by separable approximation,'' \emph{IEEE Trans. Signal
  Process.}, vol.~57, no.~7, pp. 2479--2493, Jul. 2009.

\bibitem{Wen10}
Z.~W. {Wen}, W.~{Yin}, D.~{Goldfarb}, and Y.~{Zhang}, ``A fast algorithm for
  sparse reconstruction based on shrinkage, subspace optimization, and
  continuation,'' \emph{SIAM J. Sci. Comput.}, vol.~32, no.~4, pp. 1832--1857,
  2010.

\bibitem{Liu11}
Y.~{Liu}, Y.~{Dai}, and Z.~{Luo}, ``Coordinated beamforming for miso
  interference channel: Complexity analysis and efficient algorithms,''
  \emph{IEEE Trans. Signal Process.}, vol.~59, no.~3, pp. 1142--1157, Mar.
  2011.

\bibitem{Huang15}
Y.~K. {Huang}, H.~W. {Liu}, and S.~S. {Zhou}, ``Quadratic regularization
  projected barzilai--borwein method for nonnegative matrix factorization,''
  \emph{Data min. knowl. discov.}, vol.~29, no.~6, pp. 1665--1684, 2015.

\bibitem{Tan16}
C.~H. {Tan}, S.~Q. {Ma}, Y.~H. {Dai}, and Y.~Q. {Qian}, ``Barzilai-borwein step
  size for stochastic gradient descent,'' in \emph{Adv. Neur. Inform. Process.
  Syst.}, 2016, pp. 685--693.

\bibitem{Tan10}
L.~T. {Tan}, H.~Y. {Kong}, and V.~N.~Q. {Bao}, ``Projected barzilai-borwein
  methods applied to distributed compressive spectrum sensing,'' in \emph{2010
  IEEE Symp. New Front. Dynamic Spectrum (DySPAN)}, Apr. 2010, pp. 1--7.

\bibitem{Deroo12}
F.~{Deroo}, M.~{Ulbrich}, B.~D.~O. {Anderson}, and S.~{Hirche}, ``Accelerated
  iterative distributed controller synthesis with a barzilai-borwein step
  size,'' in \emph{2012 IEEE 51st IEEE Conf. Decis. Control (CDC)}, Dec. 2012,
  pp. 4864--4870.

\bibitem{Bubeck15}
S.~{Bubeck}, ``Convex optimization: Algorithms and complexity,'' \emph{Found.
  Trends Mach. Learn.}, vol.~8, no. 3-4, pp. 231--357, Nov. 2015.

\bibitem{Horn12}
R.~A. Horn and C.~R. Johnson, \emph{Matrix analysis}.\hskip 1em plus 0.5em
  minus 0.4em\relax Cambridge university press, 2012.

\bibitem{XiADD18}
C.~{Xi}, R.~{Xin}, and U.~A. {Khan}, ``Add-opt: Accelerated distributed
  directed optimization,'' \emph{IEEE Transactions on Automatic Control},
  vol.~63, no.~5, pp. 1329--1339, May 2018.

\bibitem{Erdos59}
P.~{Erdos} and A.~{Renyi}, ``On random graphs i,'' \emph{Publ. Math. Debrecen},
  vol.~6, pp. 290--297, 1959.

\bibitem{BoydW04}
S.~{Boyd}, P.~{Diaconis}, and L.~{Xiao}, ``Fastest mixing markov chain on a
  graph,'' \emph{SIAM rev.}, vol.~46, no.~4, pp. 667--689, 2004.

\bibitem{Li17}
Z.~{Li}, W.~{Shi}, and M.~{Yan}, ``A decentralized proximal-gradient method
  with network independent step-sizes and separated convergence rates,''
  \emph{IEEE Trans. Signal Process.}, Apr. 2019.

\end{thebibliography}
\end{document}